\documentclass[final]{siamltex}
\usepackage{graphicx} %incorporate graphics in pdf, jpg, png or tif(with one f) formats
\usepackage{amsmath,amssymb}
\newcommand{\adots}{\mathinner{\mskip1mu\raise
1pt\vbox{\kern7pt\hbox{.}}\mskip2mu\raise4pt\hbox{.}
\mskip1mu\raise7pt\hbox{.}\mskip1mu}}
\title{Truncated decompositions and
filtering methods with Reflective/Anti-Reflective boundary
conditions: a comparison}
\author{C.~Tablino~Possio\thanks{Dipartimento di Matematica e Applicazioni,
Universit\`a di Milano Bicocca, via Cozzi 53, 20125 Milano, Italy
(cristina.tablinopossio@unimib.it). %
The work of the author was partially supported by
MIUR 2006017542}%
}
\begin{document}
\maketitle
%\hsize=16truecm \hoffset=-3cm \noindent
%-------------------------------------------------------------------------------------------------------
\begin{abstract}
The paper analyzes and compares some spectral filtering methods as
truncated singular/eigen-value decompositions and
Tikhonov/Re-blurring regularizations in the case of the recently
proposed Reflective \cite{NCT-SISC-1999} and Anti-Reflective
\cite{S-SISC-2003} boundary conditions.  We give numerical
evidence to the fact that spectral decompositions (SDs) provide a
good image restoration quality and this is true in particular for
the Anti-Reflective SD, despite the loss of orthogonality in the
associated transform. The related computational cost is comparable
with previously known spectral decompositions, and results
substantially lower than the singular value decomposition. The
model extension to the cross-channel blurring phenomenon of color
images is also considered and the related spectral filtering
methods are suitably adapted.\newline
\end{abstract}
\pagestyle{myheadings} \thispagestyle{plain}
\markboth{C.~Tablino~Possio}{Truncated decompositions and
filtering methods with R/AR BCs}
%-----------------------------------Section----------------------------------------------------------------
\section{INTRODUCTION} %/THE IMAGE RESTORATION PROBLEM ***
In this paper we deal with the classical image restoration problem
of blurred and noisy images in the case of a space invariant
blurring. Under such assumption the image formation process is
modelled according to the following integral equation with space
invariant kernel
\begin{equation} \label{eq:modello_continuo}
g(x)=\int h(x-\tilde x) f(\tilde x) d\tilde x + \eta(x),\ x \in
\mathbb{R}^2,
\end{equation}
where $f$ denotes the true physical object to be restored, $g$ is
the recorded blurred and noisy image, $\eta$ takes into account
unknown errors in the collected data, e.g. measurement
%or discretization
errors and noise.\newline As customary, we consider the
discretization of (\ref{eq:modello_continuo}) by means of a
standard $2D$ generalization of the rectangle quadrature formula
on an equispaced grid, ordered row-wise from the top-left corner
to the bottom-right one. Hence, we obtain the relations
\begin{equation} \label{eq:modello_discreto_base}
%g(x_{i})=\sum_{j \in \mathbb{Z}^2} H(x_{i},x_{j})f(x_{j})\Delta
%x_1\Delta x_2.
g_{i}=\sum_{j \in \mathbb{Z}^2}h_{i-j}f_{j}+\eta_i,\quad i \in
\mathbb{Z}^2,
\end{equation}
in which an infinite and a shift-invariant matrix $\widetilde
A_\infty=[h_{i-j}]_{(i,j)=((i_1,i_2),(j_1,j_2))}$, i.e., a
two-level Toeplitz matrix, is involved.\newline
%Clearly, we are considering here the
%standard Toeplitz notation where the first index denotes the
%blocks and the second the position inside the block.***\newline
%
In principle, (\ref{eq:modello_discreto_base}) presents an
infinite summation since the true image scene does not have a
finite boundary. Nevertheless, the data $g_i$ are clearly
collected only at a finite number of values, so representing only
a finite region of such an infinite scene. In addition, the
blurring operator typically shows a finite support, so that it is
completely described by a Point Spread Function (PSF) mask such as
\begin{equation} \label{eq:coefficienti_psf}
h_{PSF}
=\left[h_{i_1,i_2}\right]_{i_1=-q_1,\ldots,q_1,i_2=-q_2,\ldots,q_2}
\end{equation}
where $h_{i_1,i_2}\ge 0$ for any $i_1,i_2$ and
$\sum_{i=-q}^{q}h_{i} =1$, $i=(i_1,i_2)$, $q=(q_1,q_2)$
(normalization according to a suitable conservation law).\newline
Therefore, relations (\ref{eq:modello_discreto_base}) imply
\begin{equation} \label{eq:modello_discreto}
g_{i}= \sum_{s=-q}^{q}h_{s}f_{i-s} +\eta_i, \quad i_1=1, \ldots,
n_1, i_2=1, \ldots, n_2,
\end{equation}
where the range of collected data defines the so called  Field of
View (FOV).\newline Once again, we are assuming that all the
involved data in (\ref{eq:modello_discreto_sistema}), similarly to
(\ref{eq:modello_discreto_base}), are reshaped in a row-wise
ordering. In such a way we obtain  the linear system
\begin{equation} \label{eq:modello_discreto_sistema}
\widetilde A \tilde f =g -\eta
\end{equation}
where $\widetilde A\in \mathbb{R}^{N(n)\times N(n+2q)}$ is a
finite principal sub-matrix of $\widetilde A_\infty$, with main
diagonal containing $h_{0,0}$, $\tilde f\in \mathbb{R}^{N(n+2q)}$,
$g,\eta \in \mathbb{R}^{N(n)}$ and with $N(m)=m_1m_2$, for any
two-index $m=(m_1,m_2)$. \newline Such a reshape is considered
just to perform the theoretical analysis, since all the
de\-blur\-ring/de\-noi\-sing methods are able to deal directly
with data in matrix form. For instance, it is evident that the
blurring process in (\ref{eq:modello_discreto}) consists in a
discrete convolution between the PSF mask, after a rotation of
$180^\circ$, and the proper true image data in
\[
\widetilde{F} =\left[ f_{i_1,i_2}
\right]_{i_1=-q_1+1,\ldots,n_1+q_1,i_2=-q_2+1,\ldots,n_2+q_2}.
\]
Hereafter, with a two-index notation, we denote by $F =\left[
f_{i_1,i_2} \right]_{i_1=1,\ldots,n_1,i_2=1,\ldots,n_2}$ the true
image inside the FOV and by $G =\left[ g_{i_1,i_2}
\right]_{i_1=1,\ldots,n_1,i_2=1,\ldots,n_2}$ the recorded image.
Thus, assuming the knowledge of PSF mask  in
(\ref{eq:coefficienti_psf}) and of some statistical properties of
$\eta$, the deblurring problem is defined as to restore, as best
as possible, the true image $F$ on the basis of the recorded image
$G$. As evident from (\ref{eq:modello_discreto}), the problem is
undetermined since the number of unknowns involved in the
convolution exceeds the number of recorded data. Boundary
conditions (BCs) are introduced to artificially describe the scene
outside the FOV: the values of unknowns outside the FOV are fixed
or are defined as linear combinations of the unknowns inside the
FOV, the target being to reduce
(\ref{eq:modello_discreto_sistema}) into a square linear system
\begin{equation} \label{eq:sistema_quadrato}
A_n f = g -\eta
\end{equation}
with $A_n \in \mathbb{R}^{N(n)\times N(n)}$, $n=(n_1,n_2)$,
$N(n)=n_1 n_2$ and $f,g,\eta \in \mathbb{R}^{N(n)}$.\newline The
choice of the BCs  does not affect the global spectral behavior of
the matrix. However, it may have a valuable impact both with
respect to the accuracy of the restored image and to the
computational costs for recovering $f$ from the blurred datum,
with or without noise.
\newline Notice also that, typically, the matrix $A$ is very
ill-conditioned and there is a significant intersection between
the subspace related to small eigen/singular values and the high
frequency subspace.
%the degenerating subspace largely intersects the high frequency space.
Such a feature requires the use of suitable regularization methods
that allow to properly restore the image $F$ with controlled noise
levels \cite{EHN-2000,G-1984,H-1997,V-2002}, among which we can
cite truncated SVD, Tikhonov, and total variation
\cite{EHN-2000,H-1997,V-2002}.\newline Hereafter, we focus our
attention on special case of PSFs satisfying a strong symmetry
property, i.e., such that
\begin{equation}\label{eq:strong_simmetry_condition}
h_{|i|}=h_i \quad \textrm{for any } i=-q,\ldots,q.
\end{equation}
This assumption is fulfilled in the majority of models in real
optical applications. For instance, in most 2D astronomical
imaging with optical lens \cite{BB-1998} the model of the PSF is
circularly symmetric, and hence, strongly symmetric; in the
multi-image deconvolution of some recent interferometric
telescopes, the PSF is strongly symmetric too \cite{BB-AAS-2000}.
Moreover, in real applications when the PSF is obtained by
measurements (like a guide star in astronomy), the influence of
noise leads to a numerically nonsymmetric PSF, also when the
kernel of the PSF is strongly (or centro) symmetric. In such a
case, by employing a symmetrized version of the measured PSF,
comparable restorations are observed
\cite{HN-IP-1996,ABBCL-AA-2006}.
%In 3D out-of-focus deconvolution
%microscopy \cite{dec3d}, the PSF is again strongly symmetric.
%Indeed, the best optics yield PSFs that are symmetrical in the $Z$
%axes and circularly symmetric with respect to the $X\!\!-\!\!Y$
%plane. Other 3D applications with PSFs strongly symmetric can be
%found in \cite{Nagy3D}.
%----------------------------------------------------------------
\par
The paper is organized as follows. In Section \ref{sez:BC} we
focus on two recently proposed BCs, i.e., the Reflective
\cite{NCT-SISC-1999} and Anti-Reflective BCs \cite{S-SISC-2003}
and their relevant properties. Section \ref{sez:filtering_methods}
summarizes some classical filtering techniques as the truncated
singular/eigen-values decomposition and the Tikhonov method. The
Re-blurring method \cite{DS-IP-2005,DEMS-IP-2006} is considered in
the case of Anti-Reflective BCs and its re-interpretation in the
framework of the classical Tikhonov regularization is given. In
Section \ref{sez:cross_channel_blurring} the model is generalized
for taking into account the cross-channel blurring phenomenon and
the previous filtering methods are suitable adapted. Lastly,
Section \ref{sez:numerical_tests} deals with some computational
issues and reports several numerical tests, the aim being  to
compare the quoted filtering methods and the two type of BCs, both
in the case of gray-scale and color images. In Section
\ref{sez:conclusioni} some conclusions and remarks end the
paper.\newline
%-----------------------------------Section----------------------------------------------------------------
\section{BOUNDARY CONDITIONS} \label{sez:BC}
In this section we summarize the relevant properties of two
recently proposed type of BCs, i.e., the Reflective
\cite{NCT-SISC-1999} and Anti-Reflective BCs \cite{S-SISC-2003}.
Special attention is given to the structural and spectral
properties of the arising matrices. In fact, though the choice of
the BCs does not affect the global spectral behavior of the matrix
$A$, it can have a valuable impact with respect  both to the
accuracy of the restoration (especially close to the boundaries
where ringing effects can appear), and the computational costs for
recovering the image from the blurred one, with or without noise.
\newline Moreover, tanking into account the scale of the problem,
the regularization methods analysis can be greatly simplified
whenever a spectral (or singular value) decomposition of $A$ is
easily available. This means that the target is to obtain the best
possible approximation properties, keeping unaltered the fact that
the arising matrix shows an exploitable structure. For instance,
the use of periodic BCs enforces a circulant structure, so that
the spectral decomposition can be computed efficiently with the
fast Fourier transform (FFT) \cite{D-1979}. Despite these
computational facilities, they give rise to significant ringing
effects when a significant discontinuity is introduced into the
image. \newline Hereafter, we focus on two recently proposed
boundary conditions, that more carefully describe the scene
outside the FOV.\newline Clearly, several other methods deal with
this topic in the image processing literature, e.g. local mean
value \cite{SC-APNUM-} or extrapolation techniques (see
\cite{LB-1991} and references therein). Nevertheless, the penalty
of their good approximation properties could lie in a linear
algebra problem more difficult to cope with.
%-----------------------------------Sub-Section----------------------------------------------------------------
\subsection{REFLECTIVE BOUNDARY CONDITIONS}
In \cite{NCT-SISC-1999} Ng \emph{et al.} analyze the use of
Reflective BCs, both from model and linear algebra point of view.
The improvement with respect to Periodic BCs is due to the
preservation of the continuity of the image. In fact, the scene
outside the FOV is assumed to be a reflection of the scene inside
the FOV. For example, with a boundary at $x_1=0$ and $x_2=0$ the
reflective condition is given by $ f(\pm x_1,\pm
x_2)=f(x_1,x_2)$.\newline More precisely, along the borders, the
BCs impose %, for any $i_1=1,\ldots,n_1,\ i_2=1,\ldots, q_2$
\[
\begin{array}{rcl rcl l}
  f_{i_1,1-i_2}   &\!\!=\!\!& f_{i_1,i_2},       &f_{i_1,n_2+i_2} &\!\!=\!\!& f_{i_1,n_2+1-i_2}, &\!\!\!\! \textrm{for any } i_1=1,\ldots,n_1,\ i_2=1,\ldots, q_2\\
  f_{1-i_1,i_2}   &\!\!=\!\!& f_{i_1,i_2},       & f_{n_1+i_1,i_2} &\!\!=\!\!& f_{n_1+1-i_1,i_2}, &\!\!\!\! \textrm{for any } i_1=1,\ldots, q_1,\ i_2=1,\ldots,n_2, \\
\end{array}
\]
and, at the corners, the BCs impose for any  $i_1=1,\ldots,q_1$,
$i_2=1,\ldots,q_2$
\[
\begin{array}{rcl c rcl}
 f_{1-i_1,1-i_2}     &=& f_{i_1,i_2},       &\quad & f_{n_1+i_1,n_2+i_2} &=& f_{n_1+1-i_1,n_2+1-i_2}, \\
 f_{1-i_1,n_2+i_2}   &=& f_{i_1,n_2+1-i_2}, &\quad & f_{n_1+i_1,1-i_2}   &=& f_{n_1+1-i_1,i_2}, \\
\end{array}
\]
i.e., a double reflection, first with respect to one axis and
after with respect to the other, no matter about the
order.\newline As a consequence the rectangular matrix $\widetilde
A$ is reduced to a square Toeplitz-plus-Hankel block matrix with
Toeplitz-plus-Hankel blocks, i.e., $A_n$ shows the two-level
Toeplitz-plus-Hankel structure.
Moreover, if the blurring operator satisfies the strong symmetry
condition (\ref{eq:strong_simmetry_condition}) then the matrix
$A_n$ belongs to DCT-III matrix algebra. Therefore, its spectral
decomposition can be computed very efficiently using the fast
discrete cosine transform (DCT-III) \cite{S-1999}.\newline More in
detail, let $\mathcal{C}_n=\{A_n \in \mathbb{R}^{N(n)\times N(n)},
n=(n_1,n_2), N(n)=n_1 n_2\ | \  A_n= R_n \Lambda_n R_n^T \}$ be
the two-level DCT-III  matrix algebra, i.e., the algebra of
matrices that are simultaneously diagonalized by the orthogonal
transform
\begin{equation}\label{Rn}
R_n=R_{n_1}\otimes R_{n_2}, \quad
R_m=\left[ \sqrt{\frac{2-\delta_{t,1}}{m}}
\cos\left\{\frac{(s-1)(t-1/2)\pi}{m}\right\} \right]_{s,t=1}^m,
\end{equation}
with $\delta_{s,t}$ denoting the Kronecker symbol.\newline
Thus,
the explicit structure of the matrix is
$A_n=\mathrm{Toeplitz}(V)+\mathrm{Hankel}(\sigma(V),J\sigma(V))$,
with $V=[V_{0}\ V_{1}\ \ldots\  V_{q_1}\ 0 \ldots 0]$ and where
each $V_{i_1}$, $i_1=1,\ldots, q_1$ is the unilevel DCT-III matrix
associated to the $i_1^{th}$ row of the PSF mask, i.e.,
$V_{i_1}=\mathrm{Toeplitz}(v_{i_1})+\mathrm{Hankel}(\sigma(v_{i_1}),J\sigma(v_{i_1}))$,
with $v_{i_1}=[h_{i_1,0},\ldots , h_{i_1,q_2},0,\ldots ,0]$. Here,
we denote by $\sigma$ the shift operator such that
$\sigma(v_{i_1})=[h_{i_1,1},\ldots , h_{i_1,q_2},0,\ldots ,0]$ and
by $J$ the usual flip matrix; at the block level the same
operations are intended in block-wise sense.\newline
Beside this structural characterization, the spectral description
is completely known. In fact, let $f$ be the bivariate generating
function associated to the PSF mask (\ref{eq:coefficienti_psf}),
that is
\begin{eqnarray}
     f(x_1,x_2)&=&h_{0,0}+2\sum_{s_1=1}^{q_1}h_{s_1,0} \cos(s_1x_1)
              +2\sum_{s_2=1}^{q_2}h_{0,s_2} \cos(s_2x_2) \nonumber \\
              && \ +4\sum_{s_1=1}^{q_1}\sum_{s_2=1}^{q_2}h_{s_1,s_2}
              \cos(s_1x_1)\cos(s_2x_2), \label{eq:funzione_generatrice}
\end{eqnarray}
then the eigenvalues of the corresponding matrix  $A_n\in
\mathcal{C}_n$  are given by
\[
\lambda_{s}(A_n)= f\left(x_{s_1}^{[n_1]},x_{s_2}^{[n_2]}\right),\
s=(s_1,s_2), \quad x_r^{[m]}=\frac{(r-1)\pi}{m},
\]
where $s_1=1,\ldots,n_1$, $s_2=1,\ldots,n_2$, and where the
two-index notation highlights the tensorial structure of the
corresponding eigenvectors.\newline Lastly, notice that standard
operations like matrix-vector products, resolution of linear
systems and eigenvalues evaluations can be performed by means of
FCT-III \cite{NCT-SISC-1999} within $O(n_1n_2 \log (n_1 n_2))$
arithmetic operations (ops). For example, by multiplying by
$e_1=[1, 0, \ldots,0]^T$ both the sides of
%-----------------------------------------------------------------
%A_n= R_n D_n R_n^T
%-----------------------------------------------------------------
$R_n^T A_n =  \Lambda_n R_n^T$, it holds that
\[
[\Lambda_n]_{(i_1,i_2)}= [R_n^T (A_n e_1)]_{(i_1,i_2)} / [R_n^T
e_1]_{(i_1,i_2)}, \quad i_1=1,\ldots,n_1, i_2=1,\ldots,n_2,
\]
i.e., it is enough to consider an inverse FCT-III applied to the
first column of $A_n$, with a computational cost of $O(n_1n_2 \log
(n_1 n_2))$ ops.
%
% Risoluzione in sez 3??????????????????????
%-----------------------------------Sub-Section----------------------------------------------------------------
\subsection{ANTI-REFLECTIVE BOUNDARY CONDITIONS} \label{sez:BC-AR}
More recently, Anti-reflective boundary conditions (AR-BCs) have
been proposed in \cite{S-SISC-2003} and studied
\cite{ADS-AR1,ADS-AR2,ADNS-AR3,DEMS-IP-2006,DENPS-SPIE-2003,P-NLA-2006}.
The  improvement  is due to the fact that not only the continuity
of the image, but also of the normal derivative, are guaranteed at
the boundary. This regularity, which is not shared with Dirichlet
or periodic BCs, and only partially shared with reflective BCs,
significantly reduces typical ringing artifacts.\newline The key
idea is simply to assume that the scene outside the FOV is the
anti-reflection of the scene inside the FOV. For example, with a
boundary at $x_1=0$  the anti-reflexive condition impose
$
  f(-x_1,x_2) - f(x_1^*,x_2) =- (f(x_1,x_2) -f(x_1^*,x_2)), \
  \textrm{for any } x_2,
$ where $x_1^*$ is the center of the one-dimensional
anti-reflection, i.e.,
\[f(-x_1,x_2)  = 2f(x_1^*,x_2)-
f(x_1,x_2), \textrm{for any } x_2.
\]
In order to preserve a
tensorial structure, at the corners, a double anti-reflection,
first with respect to one axis and after with respect to the
other, is considered,
%no matter about the order.
so that the BCs impose
\[
  f(-x_1,-x_2) =
  4f(x_1^*,x_2^*)-2f(x_1^*,x_2)-2f(x_1,x_2^*)+f(x_1,x_2),
 %             &= & 2f(x_1^*,-x_2)-f(x_1,-x_2) \\
 %             &= & 2[2f(x_1^*,x_2^*)-f(x_1^*,x_2)]-[2f(x_1,x_2^*)-f(x_1,x_2)] \\
\]
where $(x_1^*,x_2^*)$ is the center of the two-dimensional
anti-reflection.\newline
More precisely, by choosing as center of the anti-reflection the
first available data, along the borders, the BCs impose
\[
\begin{array}{ll l}
\!\!   f_{1-i_1,i_2}\!\!\!=\!2f_{1,i_2}\!-\!f_{i_1+1,i_2},\
&\!\!\! f_{n_1+i_1,i_2}\!\!\!=\!
  2f_{n_1,i_2}\!-\!f_{n_1-i_1,i_2}, &\!\!\!  i_1=1,\ldots, q_1, \ \!   i_2=1,\ldots,n_2,\\
\!\!  f_{i_1,1-i_2}\!\!\!=\!2f_{i_1,1}\!-\!f_{i_1,i_2+1}, \
&\!\!\! f_{i_1,n_2+i_2}\!\!\!
=\!2f_{i_1,n_2}\!-\!f_{i_1,n_2-i_2}, &\!\!\! i_1=1,\ldots, n_1, \
\!  i_2=1,\ldots,q_2.
\end{array}
\]
At the corners, the BCs impose for any $i_1=1,\ldots,q_1$ and
$i_2=1,\ldots,q_2$,
\[
\begin{array}{rcl}
 f_{1-i_1,1-i_2}     &=& 4f_{1,1}-2f_{1,i_2+1}-2f_{i_1+1,1}+f_{i_1+1,i_2+1}, \\
 f_{1-i_1,n_2+i_2}   &=& 4f_{1,n_2}-2f_{1,n_2-i_2}-2f_{i_1+1,n_2}+f_{i_1+1,n_2-i_2}, \\
 f_{n_1+i_1,1-i_2}   &=& 4f_{n_1,1}-2f_{n_1,i_2+1}-2f_{n_1-i_1,1}+f_{n_1-i_1,i_2+1},\\
 f_{n_1+i_1,n_2+i_2} &=& 4f_{n_1,n_2}-2f_{n_1,n_1-i_2}-2f_{n_1-i_1,n_2}+f_{n_1-i_1,n_2-i_2}. \\
 \end{array}
\]
As a consequence the rectangular matrix $\widetilde A$ is reduced
to a square Toeplitz-plus-Hankel block matrix with
Toeplitz-plus-Hankel blocks, plus an additional structured low
rank matrix. \newline Moreover, under the assumption of strong
symmetry of the PSF and of a mild finite support condition
(more precisely $h_i=0$ if $|i_j|\ge n-2$, %$i=(i_1,i_2)$
for some $j\in \{1,2\}$), the resulting linear system $A_nf = g$
is such that $A_n$ belongs to the $\mathcal{AR}_n^{2D}$
commutative matrix algebra \cite{ADS-AR2}. This new algebra shares
some properties with the $\tau$ (or DST-I) algebra
\cite{BC-1983}.\newline Going inside the definition, a matrix
$A_n\in \mathcal{AR}_n^{2D}$ has the following block structure
%struttura
\[
    A_n =
    \left[
    \begin{array}{c|c|c}
    D_{0}+Z^{[1]}         & 0^T   &  0 \\
    \hline
    D_{1}+Z^{[2]}         &                                          &  0 \\
    \vdots                &                                          & \vdots \\
    D_{q_1-1}+Z^{[q_1]}   &                                          &  0 \\
    D_{q_1}               &        \  \tau(D_0,\ldots,D_{q_1}) \   &  D_{q_1}\\
    0                     &                                          &  D_{q_1-1}+Z^{[q_1]} \\
    \vdots                &                                          & \vdots \\
    0                     &                                          & D_{1}+Z^{[2]}\\
    \hline
    0                     & 0^T & D_{0}+Z^{[1]} \\
    \end{array}
    \right],
\]
where $\tau(D_0,\ldots,D_{q_1})$ is a block $\tau$ matrix with
respect to the $\mathcal{AR}^{1D}$ blocks $D_{i_1}$,
$i_1=1,\ldots, q_1$ and $Z^{[k]}= 2\sum_{t=k}^{q_1} D_{t}$  for
$k=1,\ldots,q_1$. In particular, the $\mathcal{AR}^{1D}$ block
$D_{i_1}$ is associated to $i_1^{th}$ row of the PSF, i.e.,
$h_{i_1}^{[1D]}=[h_{i_1,i_2}]_{i_2=-q_2,\ldots,q_2}$ and it is
defined as
\[
D_{i_1} %=A(h_{i_1}^{(1D)})
=
    \left[
    \begin{array}{c|c|c}
    h_{i_1,0}+z_{i_1}^{[1]}     & 0^T   &  0 \\
    \hline
    h_{i_1,1}+z_{i_1}^{[2]}     &                                               & 0 \\
    \vdots                      &                                               & \vdots \\
    h_{i_1,q_2-1}+z_{i_1}^{[q_2]} &                                             & 0 \\
    h_{i_1,q_2}                 &         \ \tau(h_{i_1,0},\ldots,h_{i_1,q_2})\ & h_{i_1,q_2}  \\
    0                           &                                               & h_{i_1,q_2-1}+z_{i_1}^{[q_2]} \\
    \vdots                      &                                               & \vdots\\
    0                           &                                               & h_{i_1,1}+z_{i_1}^{[2]}\\
    \hline
    0                           & 0^T   & h_{i_1,0}+z_{i_1}^{[1]} \\
    \end{array}
    \right],
\]
where $z_{i_1}^{[k]}= 2\sum_{t=k}^{q_2} h_{i_1,t}$  for
$k=1,\ldots,q_2$ and $\tau(h_{i_1,0},\ldots,h_{i_1,q_2})$ is the
unilevel $\tau$ matrix associated to the one-dimensional PSF
$h_{i_1}^{[1D]}$ previously defined.\newline Notice that the
rank-1 correction given by the elements $z_{i_1}^{[k]}$  pertains
to the contribution of the anti-reflection centers with respect to
the vertical borders, while the low rank correction given by the
matrices $Z_{[k]}$ pertains to the contribution of the
anti-reflection centers with respect to the horizontal
borders.\newline It is evident from the above matrix structure
that favorable computational properties are guaranteed also by
virtue of the $\tau$ structure. Therefore, firstly we recall the
relevant properties of the two-level $\tau$ algebra
\cite{BC-1983}.\newline
%-------------------------TAU----------------------------------------------
Let $\mathcal{T}_n=\{A_n \in \mathbb{R}^{N(n)\times N(n)},
n=(n_1,n_2), N(n)=n_1 n_2\ | \ A_n= Q_n \Lambda_n Q_n \}$ be the
two-level $\tau$ matrix algebra, i.e., the algebra of matrices
that are simultaneously diagonalized by the symmetric orthogonal
transform
\begin{equation}\label{Qn}
Q_n=Q_{n_1}\otimes Q_{n_2}, \quad
Q_m=\left[ \sqrt{\frac{2}{m+1}}
\sin\left\{\frac{st\pi}{m+1}\right\} \right]_{s,t=1}^m.
\end{equation}
With the same notation as the DCT-III algebra case, the explicit
structure of the matrix is two level Toeplitz-plus-Hankel. More
precisely,
$$A_n=\mathrm{Toeplitz}(V)-\mathrm{Hankel}(\sigma^2(V),J\sigma^2(V))$$
with $V=[V_{0}\ V_{1}\ \ldots\  V_{q_1}\ 0 \ldots 0]$, where each
$V_{i_1}$, $i_1=1,\ldots, q_1$ is a the unilevel $\tau$ matrix
associated to the $i_1^{th}$ row of the PSF mask, i.e.,
$V_{i_1}\!=\!\mathrm{Toeplitz}(v_{i_1})-\mathrm{Hankel}(\sigma^2(v_{i_1}),J\sigma^2(v_{i_1}))$
with $v_{i_1}=[h_{i_1,0},\ldots , h_{i_1,q_2},0,\ldots ,0]$. Here,
we denote by $\sigma^2$ the double shift operator such that
$\sigma^2(v_{i_1})=[h_{i_1,2},\ldots , h_{i_1,q_2},0,\ldots ,0]$;
at the block level the same operations are intended in block-wise
sense.\newline
Once more, the spectral characterization is completely known since
for any $A_n\in \mathcal{T}_n$ the related eigenvalues are given
by
\[
\lambda_{s}(A_n)= f\left(x_{s_1}^{[n_1]},x_{s_2}^{[n_2]}\right),
s=(s_1,s_2), \quad x_r^{[m]}=\frac{r \pi}{m+1},
\]
where $s_1=1,\ldots,n_1$, $s_2=1,\ldots,n_2$, and $f$ is the
bivariate generating function associated to the $PSF$ defined in
(\ref{eq:funzione_generatrice}).
% and the two-index notation again
%highlights the tensorial structure of the
%corresponding eigenvector.
\newline
As in the DCT-III case, standard operations like matrix-vector
products, resolution of linear systems and eigenvalues evaluations
can be performed by means of FST-I within $O(n_1n_2 \log (n_1
n_2))$ (ops). For instance,
%by multiplying
%by $e_1=[1, 0, \ldots,0]^T$ both the sides of
%-----------------------------------------------------------------
%A_n= Q_n D_n Q_n
%$Q_n A_n =  \Lambda_n Q_n$, it holds that
%-----------------------------------------------------------------
it is enough to consider a FST-I applied to the first column of
$A_n$ to obtain the eigenvalues
\[
[\Lambda_n]_{(i_1,i_2)}= [Q_n (A_n e_1)]_{(i_1,i_2)} / [Q_n
e_1]_{(i_1,i_2)}, \quad i_1=1,\ldots,n_1, i_2=1,\ldots,n_2.
\]
%-------------------------AR----------------------------------------------
Now, with respect to the $\mathcal{AR}_n^{2D}$ matrix algebra, a
complete spectral characterization is given in
\cite{ADS-AR2,ADNS-AR3}. A really useful fact
%welcome discovery
is the existence of a transform $T_n$ that simultaneously
diagonalizes all the matrices belonging to $\mathcal{AR}_n^{2D}$,
although the orthogonality property is partially lost.
\begin{theorem} {\rm \cite{ADNS-AR3}}
Any matrix $A_n\in \mathcal{AR}_n^{2D}$, $n=(n_1,n_2)$, can be
diagonalized by $T_n$, i.e.,
\[
A_n=T_n \Lambda_n \widetilde T_n,\quad \widetilde{T}_n= T_n^{-1}
\]
where $T_n=T_{n_1}\otimes T_{n_2}$, $\widetilde T_n= \widetilde
T_{n_1}\otimes \widetilde T_{n_2}$, with
\[
T_m=\left [
\begin{array}{ccc}
\alpha_m^{-1} & 0^T& 0 \\
\ \\
\alpha_m^{-1}p & Q_{m-2} & \alpha_m^{-1}Jp \\
\ \\
0 & 0^T& \alpha_m^{-1}
\end{array}
\right ]
\quad and \quad
\widetilde{T}_m=\left [
\begin{array}{ccc}
\alpha_m & 0^T & 0\\
\ \\
-Q_{m-2}p & Q_{m-2} & -Q_{m-2}Jp \\
\ \\
0 & 0^T & \alpha_m
\end{array}
\right ]
\]
The entries of the vector $p\in\mathbb{R}^{m-2}$ are defined as
$p_j=1-{j}/{(m-1)}$, $j=1,\ldots,m-2$, $J\in\mathbb{R}^{m-2\times
m-2}$ is the flip matrix, and $\alpha_m$ is a normalizing factor
chosen such that the Euclidean norm of the first and last column
of $T_m$ will be equal to $1$.
\end{theorem}
\begin{theorem} {\rm \cite{ADS-AR2}} \label{teo:teo_autovalori}
Let $A_n\in \mathcal{AR}_n^{2D}$, $n=(n_1,n_2)$, the matrix
related to the PSF $h_{PSF}=[h_{i_1,i_2}]_{i_1=-q_1,\ldots,q_1,
i_2=-q_2,\ldots,q_2}$. Then, the eigenvalues of $A_n$ are given by
\begin{itemize}
\item $1$ with algebraic multiplicity $4$,
\item the $n_2-2$ eigenvalues of the  unilevel $\tau$ matrix
related to the one-dimensional PSF
$h^{\{r\}}=[\sum_{i_1=-q_1}^{q_1}h_{i_1,-q_2}, \ldots,
\sum_{i_1=-q_1}^{q_1}h_{i_1,q_2}]$,
%(concentrazione rispetto alle righe della PSF $h$ - relativa ai bordi orizzontali)
each one with algebraic multiplicity $2$,
\item the $n_1-2$ eigenvalues of the  unilevel $\tau$ matrix
related to the one-dimensional PSF $h^{\{c\}}=
[\sum_{i_2=-q_2}^{q_2}h_{-q_1,i_2}, \ldots,
\sum_{i_2=-q_2}^{q_2}h_{q_1,i_2}]$,
%(concentrazione rispetto alle colonne della PSF $h$ - relativa ai bordi verticali)
each one with algebraic multiplicity $2$,
\item the $(n_1-2)(n_2-2)$ eigenvalues of the  two-level $\tau$
matrix related to the two-dimensional PSF  $h_{PSF}$.
%  each one with algebraic multiplicity $1$.
\end{itemize}
\end{theorem}
Notice that the three sets  of multiple eigenvalues are exactly
related to the type of low rank correction imposed by the BCs
through the centers of the anti-reflections. More in detail, the
eigenvalues of $\tau_{n_2-2}(h^{\{r\}})$  and of
$\tau_{n_1-2}(h^{\{c\}})$ take into account the condensed PSF
information considered along the horizontal and vertical borders
respectively, while the eigenvalue equal to $1$ takes into account
the condensed information of the whole PSF at the four
corners.\newline In addition, it is worth noticing that the
spectral characterization can be completely described in terms of
the generating function associated to the $PSF$ defined in
(\ref{eq:funzione_generatrice}), simply by extending to $0$ the
standard $\tau$ evaluation grid, i.e., it holds
\[ \lambda_{s}(A_n)=
f\left(x_{s_1}^{[n_1]},x_{s_2}^{[n_2]}\right), s=(s_1,s_2),
s_j=0,\ldots,n_j, \quad x_r^{[m]}=\frac{r \pi}{m+1},
\]
where the $0-$index refers to the first/last columns of the matrix
$T_m$ \cite{ADS-AR2}.\newline See \cite{ADS-AR1,ADNS-AR3} for some
algorithms related to standard operations like matrix-vector
products, resolution of linear systems and eigenvalues evaluations
with a computational cost of $O(n_1n_2\log (n_1 n_2))$
ops.\newline It is worthwhile stressing that the computational
cost of the inverse transform is comparable with that of the
direct transform and, at least at first sight, the very true
penalty is the loss of orthogonality due to the first/last column
of the matrix $T_m$. \newline
%-----------------------------------Section----------------------------------------------------------------
\section{FILTERING METHODS} \label{sez:filtering_methods}
Owing to the ill-conditioning, the standard solution $f=A_n^{-1}g$
is not physically meaningful since it is completely corrupted by
the noise propagation from data to solution, i.e., by the so
called inverted noise. For this reason, restoration methods look
for an approximate solution with controlled noise levels: widely
considered regularization methods are obtained through spectral
filtering \cite{H-1997,HNOL-2006}. Hereafter, we consider the
truncated Singular Values Decompositions (SVDs) (or Spectral
Decompositions (SDs)) and the Tikhonov (or Re-blurring)
regularization method.
%-----------------------------------sub-section----------------------------------------------------------------
\subsection{TRUNCATED SVDs AND TRUNCATED SDs} \label{sez:filtering_methods_TSVD_TSD}
The Singular Values Decomposition (SVD) highlights a standard
perspective for dealing with the inverted noise. More precisely,
if
\[
A_n=U_n\Sigma_n V_n^T\in \mathbb{R}^{N(n)\times N(n)}
\]
is the SVD of $A_n$, i.e., $U_n$ and $V_n$ are orthogonal matrices
and $\Sigma_n$ is a diagonal matrix with entries $\sigma_1 \ge
\sigma_2 \ge \ldots \sigma_{N(n)} \ge 0$, then the solution of the
linear system $A_nf=g$ can be written as
\[
f=\sum_{k=1}^{N(n)} \left ( \frac{u_k^T g}{\sigma_k}\right) v_k,
\]
where $u_k$ and $v_k$  denote the $k^{th}$ column of the matrix
$U_n$ and $V_n$, respectively.\newline With regard to the image
restoration problem, the idea is to consider a sharp filter, i.e.,
to take in the summation only the terms corresponding to singular
values greater than a certain threshold value $\delta$, so damping
the effects caused by division by the small singular values.
Therefore, by setting the filter factors as
\[
\phi_k =\left \{
\begin{array}{lcl}
1, & \ & \mathrm{if\ } \sigma_k \ge \delta,  \\
0, & \ & \mathrm{otherwise},
\end{array}
\right.
\]
the filtered solution is defined as
\[
f_{\mathrm{filt}}=\sum_{k=1}^{N(n)} \left ( \phi_k\frac{u_k^T
g}{\sigma_k}\right) v_k =\sum_{k \in I_\delta} \left (
\phi_k\frac{u_k^T g}{\sigma_k}\right) v_k, \quad I_\delta=\{k\ | \
\sigma_k \ge \delta\}.
\]
Due to scale of the problem, the SVD of the matrix $A_n$ is in
general an expensive computational task (and not negligible also
in the case of a separable PSF). Thus, an ``a priori'' known
spectral decomposition, whenever available, can give rise to a
valuable simplification. More precisely, let
\[
A_n=V_n\Lambda_n \widetilde V_n \in \mathbb{R}^{N(n)\times
N(n)},\quad \widetilde V_n= V_n^{-1}
\]
be a spectral decomposition of $A_n$, then the filtered solution
is defined as
\[
f_{\mathrm{filt}}=\sum_{k=1}^{N(n)}
\left ( \phi_k\frac{\tilde v_k g}{\lambda_k}\right)
v_k =\sum_{k \in I_\delta}
\left ( \phi_k\frac{\tilde v_k g}{\lambda_k}\right)
v_k, \quad I_\delta=\{k\ | \ |\lambda_k(A)| \ge \delta\},
\]
where $v_k$ and $\tilde v_k$  denote the $k^{th}$ column of $V_n$
and the $k^{th}$ row of $\widetilde V_n$, respectively, and where
$\phi_k=1$ if $k\in I_\delta$, $0$ otherwise.\newline
%
%-----------------------------------sub-section----------------------------------------------------------------
\subsection{TIKHONOV AND RE-BLURRING REGULARIZATIONS}
In the classical Tikhonov regularization method, the image
filtering is obtained by looking for the solution of the following
minimization problem
\begin{equation} \label{eq:min_tikhnonov}
\min_{f} \|A_n f-g \|_2^2+\mu \| D_nf\|_2^2,
\end{equation}
where $\mu>0$ is the regularization parameter and $D_n$ is a
carefully chosen matrix (typically $D_n=I_n$ or represents the
discretization of a differential operator, properly adapted with
respect to the chosen BCs).\newline The target is to minimize the
Euclidean norm of the residual $\|A_n f-g \|_2$ without explosions
with respect to the quantity $\|D_n x\|_2$. As well know,
(\ref{eq:min_tikhnonov}) is equivalent to the solution to the
damped least square problem
\begin{equation} \label{eq:dampedNEs}
(A^T_n A_n + \mu D^T_n D_n) f =A^T_n g.
\end{equation}
In addition, the regularization Tikhonov method can be
reinterpreted in the framework of classical spectral filtering
method. For instance, in the case of $D_n=I_n$, by making use of
the SVD of $A_n=U_n\Sigma_n V^T_n$, the solution of
(\ref{eq:dampedNEs}) can be rewritten as
\[
f_{\mathrm{filt}}=V_n \Phi_n \Sigma^{-1}_n U^T_n g,
\]
where $\Phi_n=\mathrm{diag}(\phi_k)$ with
$\phi_k={\sigma_k^2}/{(\sigma_k^2+\mu)}$, $k=1,\ldots,
N(n)$.\newline
%-----------------------------------sub-section----------------------------------------------------------------
%Motivazioni reblurring\\
%Analisi reblurring\\
A severe drawback in adopting the Tikhonov regularization approach
in the case of  $A_n\in \mathcal{AR}_n^{2D}$ is due to the fact
that $A_n^T \notin \mathcal{AR}_n^{2D}$, so that all the favorable
computational properties are substantially spoiled. An alternative
approach, named Re-blurring, has been proposed in
\cite{DS-IP-2005,DEMS-IP-2006}: the proposal is to replace $A^T_n$
by $A^{'}_n$ in (\ref{eq:dampedNEs}), where $A^{'}_n$ is the
blurring matrix related to the current BCs with  a PSF rotated by
$180^{\circ}$. This approach is completely equivalent to
(\ref{eq:dampedNEs}) in the case of Dirichlet and Periodic BCs,
while the novelty concerns both Reflective BCs and Anti-Reflective
BCs, where in general $A^{'}_n \ne A^T_n$. The authors show that
the Re-blurring with anti-reflective BCs is computationally
convenient and leads to a larger reduction of the ringing effects
arising in classical deblurring schemes. From the modelling point
of view, the authors motivation relies upon the fact that
Re-blurring smoothes the noise in the right hand side of the
system, in the same manner as this happens in the case of
Dirichlet, Periodic and Reflective BCs. \newline Hereafter, we
consider an explanation of the observed approximation results. As
previously claimed, we focus our attention on the case of a
strongly symmetric PSF, so that the matrix $A^{'}_n$ equals the
matrix $A_n$. Moreover, also in this case it is evident that the
linear system
\begin{equation} \label{eq:damped_reblurring}
(A^2_n + \mu D^2_n) f =A_ng.
\end{equation}
is not equivalent to a minimization problem, again because the
matrix $A\in \mathcal{AR}_n^{2D}$ is not symmetric. Nevertheless,
the symmetrization of (\ref{eq:damped_reblurring}) can be
performed by diagonalization, so obtaining
\begin{equation} \label{eq:damped_reblurring_diag}
(\Lambda_{A,n}^2 + \mu \Lambda_{D,n}^2) \hat f =\Lambda_{A,n}
 \hat g,
\end{equation}
where $\hat f=\widetilde T_n f$ and $\hat g=\widetilde T_n g$. In
such a way (\ref{eq:damped_reblurring_diag}) is again equivalent
to the minimization problem
\begin{equation} \label{eq:min_reblurring_1}
\min_{f} \|\Lambda_{A,n} \widetilde T_n f- \widetilde T_n g
\|_2^2+\mu \| \Lambda_{D,n} \widetilde T_n f\|_2^2,
\end{equation}
or equivalently, again by making use of the diagonalization
result, to
\begin{equation} \label{eq:min_reblurring_2}
\min_{f} \| \widetilde T_n( A_n  f- g) \|_2^2+\mu \| \widetilde
T_n D_n f\|_2^2.
\end{equation}
Clearly, the last formulation in (\ref{eq:min_reblurring_2}) is
the most natural and it allows to claim that the Re-blurring
method can be interpreted as a standard Tikhonov regularization
method in the space transformed  by means of $\widetilde T_n$.
\newline Recalling that $\widetilde T_n$ is not an orthogonal
transformation, the goal becomes to compare $\| \widetilde T_n
f\|_2$ and $\|f\|_2$, that is to bound $\| \widetilde T_n\|_2=\|
\widetilde T_{n_1}\|_2\| \widetilde T_{n_1}\|_2$, being $\|
\widetilde T_n f\|_2 \le \| \widetilde T_n\|_2 \| f\|_2$. \newline
%%%%%%%%%%%%%%%%%%%%%%%%%%%%%%%%%%%%%%%%%%%%%%%%%%%%%%%%%%%%%%%%%%%
A quite sharp estimate of such a norm can be found by exploiting
the structure of the unilevel matrix $\widetilde T_m \in
\mathbb{R}^{m \times m}$. Let $\breve{f}=[f_2,\ldots,f_{m-1}]$, it
holds that
\begin{eqnarray*}
\|\widetilde T_m f\|^2_2
&=& \alpha_m^2 f_1^2 +\|Q_{m-2}(-f_1p+\breve{f}-f_n J p)\|_2^2
+\alpha_m^2 f_m^2\\
&=& \alpha_m^2 (f_1^2 +f_m^2)+\|-f_1p+\breve{f}-f_n J p\|_2^2\\
%
%&\le& \alpha_m^2 (f_1^2 +f_m^2)+(\|(f_1p\|_2+\|\breve{f}\|_2+\|f_n J p)\|_2)^2\\
%
&\le& \alpha_m^2 (f_1^2 +f_m^2)+(\|\breve{f}\|_2+(|f_1|+|f_n|)\|p\|_2)^2\\
&\le& \alpha_m^2 (f_1^2 +f_m^2)+\|\breve{f}\|_2^2
+ 3\|p\|_2^2 \|f\|_2^2 + 4\|p\|_2 \|f\|_2^2\\
%
%&=& \|f\|_2^2+(f_1^2 +f_m^2)\|p\|_2^2+\|\breve{f}\|_2^2 + 3\|p\|_2^2 \|f\|_2^2 + 4\|p\|_2 \|f\|_2^2,
%
&\le&(1+2 \|p\|_2)^2 \|f\|_2^2,
\end{eqnarray*}
being $\alpha_m^2=1 + \|p\|_2^2$.  Since, by definition,
$\|p\|_2^2\simeq m$, we have
\begin{equation} \label{eq:noma2_Ttilde}
 \|\widetilde T_m \|_2 \le 1+2 \|p\|_2 \simeq 2\sqrt{m}.
\end{equation}
Notice that the bound given in (\ref{eq:noma2_Ttilde}) is quite
sharp, since for instance $\|\widetilde T_m e_1\|^2_2$ equals $1 +
2\|p \|_2^2$.\newline
\section{CROSS-CHANNEL BLURRING} \label{sez:cross_channel_blurring}
Hereafter, we extend the analysis of the deblurring problem to the
case of color images digitalized, for instance, according to the
standard RGB system. Several techniques can be used for recording
color images, but the main problem concerns the fact that light
from one color channel can end up on a pixel assigned to another
color.
%For example, in the single-CCD approach, this occurs by
%misalignment of the Bayer filter.
The consequence of this phenomenon is called cross-channel
blurring among the three channels of the image and it sums up to
the previously analyzed blurring  of each one of the three colors,
named within-channel blurring.\newline By assuming that the
cross-channel blurring takes place after the
within-channel  blurring %(the optical blurring)
of the image, that it is spatially invariant and by assuming  that
the same within-channel blurring occurs in all the three color
channels, the problem can be modelled \cite{HNOL-2006} as
\begin{equation} \label{eq:sistema_quadrato_color}
(A_{\mathrm{color}}\otimes A_n) f = g -\eta
\end{equation}
with $A_n \in \mathbb{R}^{N(n)\times N(n)}$, $n=(n_1,n_2)$,
$N(n)=n_1 n_2$, and
\[
A_{\mathrm{color}}=\left [
\begin{array}{ccc}
  a_{rr} & a_{rg} & a_{rb} \\
  a_{gr} & a_{gg} & a_{gb} \\
  a_{br} & a_{bg} & a_{bb} \\
\end{array}
\right ].
\]
The row-entries denote the amount of within-channel blurring
pertaining to each color channel; a normalized conservation law
prescribes that $A_{\mathrm{color}} e= e$, $e=[1\ 1\ 1]^T$.
Lastly, the vectors $f,g,\eta \in \mathbb{R}^{3N(n)}$ are assumed
to collect the three color channels in the RGB order.\newline
Clearly, if $A_{\mathrm{color}} = I_3$, i.e., the blurring is only
of within-channel type, the problem is simply decoupled into three
independent gray-scale deblurring problems.\newline In the general
case, taking into account the tensorial structure of the whole
blurring matrix $A_{\mathrm{color}}\otimes A_n$ is evident that
the truncated SVDs and SDs can be formulated as the natural
extension of those considered in the within-blurring case. Notice
that in the case of SDs, we will consider a SVD for the matrix
$A_{\mathrm{color}}$, since it naturally assures an orthogonal
decomposition, no matter about the specific matrix, while its
computational cost is negligible with respect to the scale of the
problem. In addition, we tune the filtering strategy with respect
the spectral information given only by the matrix $A_n$, i.e., for
any fixed $\sigma_k$ (or $\lambda_k$) we simultaneously sum, or
discard, the three contribution on $f$ related to the three
singular values of $A_{\mathrm{color}}$.\newline With respect to
the Tikhonov regularization method, the approach is a bit more
involved. Under the assumption $A_n=A_n^T=V_n \Lambda_n \widetilde
V_n$, the damped least square problem
\[
[(A_{\mathrm{color}}\otimes A_n)^T(A_{\mathrm{color}}\otimes
A_n)+\mu I_{3n}] f = (A_{\mathrm{color}}\otimes A_n)^T g
\]
can be rewritten as
\begin{equation}
[(A_{\mathrm{color}}^T A_{\mathrm{color}}) \otimes V_n \Lambda_n^2
\widetilde V_n+\mu (I_{3}\otimes I_{n})] f = (A_{\mathrm{color}}
\otimes V_n \Lambda_n \widetilde V_n)^T g. \label{eq:damped_color}
\end{equation}
Thus, by setting $S_{3n}=I_3 \otimes \widetilde V_n$, $\hat
f=S_{3n} f$,  $\hat g=S_{3n} g$,
%and $S_{3n}^{-1}=I_3 \otimes V_n$
(\ref{eq:damped_color}) can be transformed in
\[
S_{3n}[(A_{\mathrm{color}}^T A_{\mathrm{color}}) \otimes V_n
\Lambda_n^2 \widetilde V_n+\mu (I_{3}\otimes I_{n})]
S_{3n}^{-1}\hat f = S_{3n} (A_{\mathrm{color}}\otimes V_n
\Lambda_n \widetilde V_n)^T S_{3n}^{-1}\hat g,
\]
so obtaining the linear system
\[
[(A_{\mathrm{color}}^T A_{\mathrm{color}}) \otimes \Lambda_n^2
+\mu (I_{3}\otimes I_{n})] \hat f = (A_{\mathrm{color}}^T \otimes
\Lambda_n )\hat g,
%
%\label{eq:damped_color_trasformato}
\]
that can easily be decoupled into $n_1n_2$ linear systems of
dimension $3$.\newline Clearly, in the case of any matrix $A_n\in
\mathcal{C}_n$, all these manipulations can be performed
 by means of an orthogonal transformation $S_{3n}$.
Notice also that the computational cost is always $O(n_1 n_2 \log
n_1 n_2)$ ops.\newline
With respect to $A_n=T_n \Lambda_n
\widetilde T_n \in \mathcal{AR}_n^{2D}$, we can consider the same
strategy by referring to the Re-blurring regularization method.
More precisely, the linear system
\[
[(A_{\mathrm{color}}^T A_{\mathrm{color}}) \otimes A_n^2 +\mu
(I_{3}\otimes I_{n})] f = (A_{\mathrm{color}}^T \otimes A_n) g
%\label{eq:damped_color_AR}
\]
can be transformed in
\[
[(A_{\mathrm{color}}^T A_{\mathrm{color}}) \otimes \Lambda_n^2
+\mu (I_{3}\otimes I_{n})] \hat f = (A_{\mathrm{color}}^T \otimes
\Lambda_n )\hat g.
%
%\label{eq:damped_color_trasformato_AR}
\]
Though the transformation $S_{3n}=I_3 \otimes \widetilde T_n$ is
not orthogonal as in the Reflective case, the obtained restored
image are fully comparable with the previous ones and the
computational cost is still $O(n_1 n_2 \log n_1 n_2))$ ops.
% Some numerical test in the case of cross-channel
%blurring are reported in Section
%\ref{sez:cross_channel_blurring_numerical}.
%-----------------------------------Section----------------------------------------------------------------
\section{NUMERICAL TESTS} \label{sez:numerical_tests}
%-----------------------------------sub-section----------------------------------------------------------------
\subsection{SOME COMPUTATIONAL ISSUES} \label{sez:computational_issues}
Before analyzing the image re\-sto\-ra\-tion results, we discuss
how the methods can work without reshaping the involved data. In
fact, the tensorial structure of the matrices, obtained by
considering Reflective and Anti-Reflective BCs, can be exploited
in depth, so that the algorithms can deal directly, and more
naturally, with the data collected in matrix form.
% as more natural in image processing.
Hereafter, we consider a two-index notation in the sense of the
previously adopted row-wise ordering.\newline In the SD case
considered in Section \ref{sez:filtering_methods_TSVD_TSD}, since
$\tilde v_k= \tilde v_{k_1}^{[n_1]} \otimes \tilde
v_{k_2}^{[n_2]}$ is represented in matrix form as $(\tilde
v_{k_1}^{[n_1]})^T \tilde v_{k_2}^{[n_2]}$, the required scalar
product can be computed as
\[
\tilde v_k g =\left [\left (\tilde v_{k_1}^{[n_1]}\right)^T \tilde
v_{k_2}^{[n_2]}\right] \odot G,
\]
where $\odot$ denotes the summation of all the involved terms
after a element-wise product.  Clearly, $v_k= v_{k_1}^{[n_1]}
\otimes v_{k_2}^{[n_2]}$ is represented in matrix form as
$v_{k_1}^{[n_1]} (\tilde v_{k_2}^{[n_2]})^T$. In a similar manner,
in the case of the SVD of $A_n$ with separable PSF $h=h_1\otimes
h_2$, we can represent $v_k= v_{k_1}^{[n_1]} \otimes
v_{k_2}^{[n_2]}$ in matrix form as $v_{k_1}^{[n_1]}
(v_{k_2}^{[n_2]})^T$ and $u_k^T=(u_{k_1}^{[n_1]} \otimes
u_{k_2}^{[n_2]})^T$ as $u_{k_1}^{[n_1]}
(u_{k_2}^{[n_2]})^T$.\newline The eigenvalues required for the SD
can be stored into a matrix $\Lambda^*\in \mathbb{R}^{n_1\times
n_2}$. In the case of $A_n\in \mathcal{C}_n$ this matrix can be
evaluated as
\[
\Lambda^*= \left( \widetilde V_{n_2} A^* \widetilde V_{n_1}^T\ ./
\ \widetilde V_{n_2} E_1^* \widetilde V_{n_1}^T \right)^T
\]
where $A^*\in \mathbb{R}^{n_2\times n_1}$ denotes the first column
of $A_n$ and $E_1^*$ the first canonical basis vector, reshaped as
matrices in column-wise order. In addition, the two-level direct
and inverse transform $y=V_nx$ and $y=\widetilde V_n x$ can be
directly evaluated on a matrix data as
\[
Y=V_{n_1}X V_{n_2}^T = (V_{n_2}(V_{n_1}X)^T)^T\quad \mathrm{and}
\quad Y=\widetilde V_{n_1}X \widetilde V_{n_2} = (\widetilde
V_{n_2}(\widetilde V_{n_1}X)^T)^T
\]
by referring to the corresponding unilevel transforms.\newline
%-----------------------------------------------------------------------
% VER ordine PSF separabile
%-----------------------------------------------------------------------
In the same way, the eigenvalues required in the case of $A_n\in
\mathcal{AR}_n^{2D}$ can be suitably stored as
\[
\Lambda^* = \left [
\begin{array}{c|c|c}
1 &  \ \Lambda^*(\tau_{n_2-2}(h^{r})) \ & 1\\
\hline
& &\\
& &\\
\Lambda^*(\tau_{n_1-2}(h^{c})) &  \Lambda^*(\tau_{n-2}(h)) &  \Lambda^*(\tau_{n_1-2}(h^{c}))\\
& &\\
& &\\
\hline
1 & \Lambda^*(\tau_{n_2-2}(h^{r}))  & 1\\
\end{array}
\right ] \in \mathbb{R}^{n_1\times n_2},
\]
with reference to the notations of Theorem
\ref{teo:teo_autovalori}, where the eigenvalues of the unilevel
and two-level $\tau$ matrices are evaluated as outlined in Section
\ref{sez:BC-AR}.\newline
%-----------------------------------------------------------------------
% CASO COSENI
%With respect to the required eigenvalues, we can consider the
%matrix
%\[
%\Lambda^*= ( C_{n_2}^T A^* C_{n_1} ./C_{n_2}^T E_1^* C_{n_1})^T
%\in \mathbb{R}^{n_1\times n_2}
%\]
%where $A^*\in \mathbb{R}^{n_2\times n_1}$ denotes the first column
%of $A$, reshaped as a matrix in column-wise order.\\
%Moreover, the two-level transform $y=C_Nx$ can be directly
%evaluated on a matrix data as
%\[
%Y=C_{n_1}X C_{n_2}^T = (C_{n_2}(C_{n_1}X)^T)^T
%\]
%by referring to the corresponding unilevel transforms.\\
%In the same way,  $y=C_N^Tx=(C_{n_1}^T\otimes C_{n_2}^T)x$  can be
%directly evaluated as
%\[
%Y=C_{n_1}^TX C_{n_2} = (C_{n_2}^T(C_{n_1}^TX)^T)^T
%\]
%by referring to the unilevel inverse transforms.\\
%
%
%-----------------------------------------------------------------------
Lastly, the linear systems obtained, for any fixed $\mu$, in the
case of Tikhonov and Re-blurring regularization methods can be
solved  with reference to the matrix $\Phi_n$ of the corresponding
filter factors by applying the Reflective and Anti-Reflective
transforms with a computational cost $O(n_1n_2\log n_1 n_2)$ ops.
%-----------------------------------------------------------------------
%In a similar manner, in the case of the SVD, we have that $u_k^T=
%(u_{k_1}^{[n_1]} \otimes u_{k_2}^{[n_2]})^T$ is represented in
%matrix form as $(u_{k_1}^{[n_1]} (u_{k_2}^{[n_2]})^T$, so that the
%scalar product can be computed as
%
%\[
%u_k^T g =\left [ u_{k_1}^{[n_1]} \left
%(u_{k_2}^{[n_2]}\right)^T\right] \odot G,
%\]
%
%where $\odot$ denotes the summation of all the involved terms
%after a element-wise product. \\
%Clearly, $v_k= v_{k_1}^{[n_1]} \otimes v_{k_2}^{[n_2]}$ is
%represented in matrix
%form as $v_{k_1}^{[n_1]} (v_{k_2}^{[n_2]})^T$.\\
%-----------------------------------sub-section----------------------------------------------------------------
\subsection{TRUNCATED DECOMPOSITIONS}
In this section we  compare the effectiveness of truncated
spectral decompositions (SDs)  with respect to the standard
truncated SVDs both in the case of  Reflective and Anti-Reflective
BCs. Due to scale of the problem, the SVD of the matrix $A_n$ is
in general an expensive computational task (and not negligible
also in the case of a separable PSF). Thus, a spectral
decomposition, whenever available as in these cases, leads to a
valuable simplification.\newline Firstly, we consider the case of
the separable PSF caused by atmospheric turbulence
\[
h_{i_1,i_2}= \frac{1}{2\pi \sigma_{i_1} \sigma_{i_2}} \exp \left (
-\frac{1}{2} \left (\frac{i_1}{\sigma_{i_1}}\right)^2
-\frac{1}{2} \left (\frac{i_2}{\sigma_{i_2}}\right)^2
\right),
\]
where $\sigma_{i_1}$ and $\sigma_{i_2}$ determine the width of the
PSF itself. Since the Gaussian function decays exponentially away
from its center, it is customary to truncate the  values in the
PSF mask after an assigned decay $|i_1|,|i_2|\le l$. It is evident
from the quoted definition that the Gaussian PSF satisfies the
strong symmetry condition (\ref{eq:strong_simmetry_condition}).
Another example of strongly symmetric PSF is given by the PSF
representing the out-of-focus blur
\[
h_{i_1,i_2}= \left \{
\begin{array}{lcl}
\frac{1}{\pi r^2}, &\quad &\textrm{if }
i_1^2+i_2^2 \le r^2,\\
0, && \mathrm{otherwise},
\end{array}
\right.
\]
where $r$ is the radius of the PSF.\newline
%%%%%%%%%%%%%%%%%%%%%%%%%%%%%%%%%%%%%%%%%%%%%%%%%%%%%%%%%%%%%%%%%%%%%%%%%%%%%%%%%%%%%
% NOISE
%randn *********
% distribuzione normale centrata in 0, con varianza e deviazione standard pari a 1.
% eta_r=ratio*norm(G,'fro')/norm(noise,'fro');
% Gnoised=G+eta*noise;
%%%%%%%%%%%%%%%%%%%%%%%%%%%%%%%%%%%%%%%%%%%%%%%%%%%%%%%%%%%%%%%%%%%%%%%%%%%%%%%%%%%%%
In the reported numerical tests, the blurred image $g$ has been
perturbed by adding a Gaussian noise contribution $\eta=\eta_n
\nu$ with $\nu$ fixed noise vector,  $\eta_n=\rho
\|g\|_2/\|\nu\|_2$, and $\rho$ assigned value. In such a way the
Signal Noise Ratio (SNR) \cite{BB-1998} is given by
\[
SNR = 20 \log_{10} \frac{\|g\|_2}{\|\eta\|_2} = 20 \log_{10}
\rho^{-1}\ (dB).
\]
%
%%%%%%%%%%%%%%%%%%%%%%%%%%%%%%%%%%%%%%%%%%%%%%%%%%%%%%%%%%%%%%%%%%%%%%%%%%%%%%%%%%%%%
\subsubsection{GRAY-SCALE IMAGES}
%%%%%%%%%%%%%%%%%%%%%%%%%%%%%%%%%%%%%%%%%%%%%%%%%%%%%%%%%%%%%%%%%%%%%%%%%%%%%%%%%%%%%
\begin{figure}[t!]
   \centering
   \includegraphics[width=\textwidth]{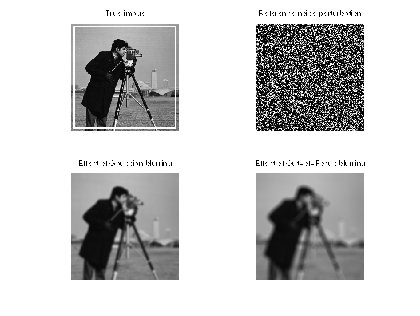}
   \caption{True image (FOV is delimited by a white frame), reference noise
perturbation,  blurred image with the Gaussian PSF with support
$15\times 15$ and $\sigma_{i_1}=\sigma_{i_2}=2$, and blurred image
with the Out-of-Focus PSF with support $15\times 15$.}
   \label{fig:cameraman_psfgauss_d15s2}
\end{figure}
\begin{figure}
   \centering
   $\rho=0.01$ \\
   \includegraphics[width=\textwidth]{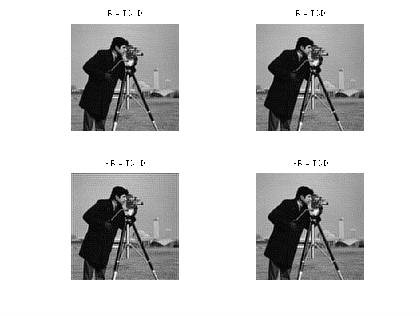}
%
%   \caption{Optimal restorations of truncated SVD and SD in the
%   case of Reflective and Anti-Reflective BCs
%   with reference to Fig. \ref{fig:cameraman_psfgauss_d15s2}.% - $\rho=0.01$.
%   }
%   \label{fig:cameraman_psfgauss_d15s2ratio001}
%\end{figure}
%%%%%%%%%%%%%%%%%%%%%%%%%%%%%%%%%%%%%%%%%%%%%%%%%%%%%%%%%%%%%%%%%%%%%%%%%%%%%%%%%%%%%
%-------------------------------------------------------------------------------------------------------------
%%%%%%%%%%%%%%%%%%%%%%%%%%%%%%%%%%%%%%%%%%%%%%%%%%%%%%%%%%%%%%%%%%%%%%%%%%%%%%%%%%%%%
%\begin{figure}
%   \centering
   $\rho=0.05$\\
   \includegraphics[width=\textwidth]{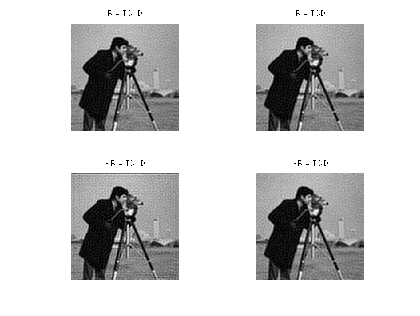}
   \caption{Optimal restorations of truncated SVD and SD in the
   case of Reflective and Anti-Reflective BCs
   with reference to Figure \ref{fig:cameraman_psfgauss_d15s2} (Gaussian blur $\sigma_{i_1}=\sigma_{i_2}=2$).% - $\rho=0.05$
   }
   \label{fig:cameraman_psfgauss_d15s2ratio001-005}
%   \label{fig:cameraman_psfgauss_d15s2ratio005}
\end{figure}

\begin{figure}
   \centering
   $\rho=0.01$ \\
   \includegraphics[width=\textwidth,height=6cm]{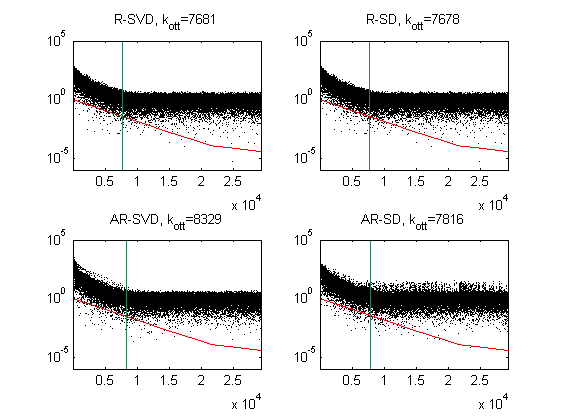}\\
%   \includegraphics[width=\textwidth, height=6cm]{immagini_Ar/fig4_picard_plot_n.eps}
%
%   \caption{Picard plot of truncated SVD and SD in the
%   case of Reflective and Anti-Reflective BCs
%   with reference to Figure \ref{fig:cameraman_psfgauss_d15s2} - $\rho=0.01$.}
%   \label{fig:cameraman_psfgauss_d15s2ratio001_picard_plot}
%\end{figure}
%%%%%%%%%%%%%%%%%%%%%%%%%%%%%%%%%%%%%%%%%%%%%%%%%%%%%%%%%%%%%%%%%%%%%%%%%%%%%%%%%%%%%
%\begin{figure}
   \centering
   $\rho=0.05$ \\
   \includegraphics[width=\textwidth, height=6cm]{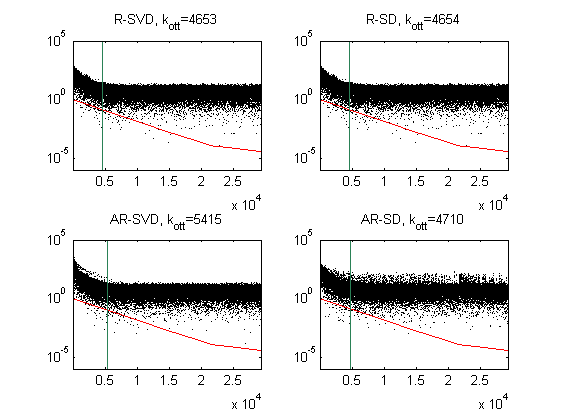}
   \caption{Picard plot of truncated SVD and SD in the
   case of Reflective and Anti-Reflective BCs
   with reference to Figure \ref{fig:cameraman_psfgauss_d15s2} (Gaussian blur $\sigma_{i_1}=\sigma_{i_2}=2$).% - $\rho=0.05$.
   }
   \label{fig:cameraman_psfgauss_d15s2ratio001-005_picard_plot}
\end{figure}

 \ \newline
%%%%%%%%%%%%%%%%%%%%%%%%%%%%%%%%%%%%%%%%%%%%%%%%%%%%%%%%%%%%%%%%%%%%%%%%%%%%%%%%%%%%%
In Figure \ref{fig:cameraman_psfgauss_d15s2} we report the
template true image (the FOV is delimited by a white frame),
together with the blurred image with the Gaussian PSF with support
$15\times 15$ and $\sigma_{i_1}=\sigma_{i_2}=2$ and the reference
perturbation $\nu$, reshaped in matrix form.\newline We consider
the optimal image restoration with respect to the relative
restoration error (RRE), i.e., $\|f_{\mathrm{filt}} -
f_{\mathrm{true}}\|_ 2/\| f_{\mathrm{true}} \|_2$, where
$f_{\mathrm{filt}}$ is the computed approximation of the true
image $f_{\mathrm{true}}$ by considering spectral filtering. More
in detail, the RRE is analyzed by progressively adding a new basis
element at a time, according to the non-decreasing order of the
singular/eigen-values (the eigenvalues are ordered with respect to
their absolute value).\newline In the case of SDs (or SVDs related
to a separable PSF) this can be done as described in Section
\ref{sez:computational_issues} and, beside the preliminary cost
related to the decomposition computation, the addition of a new
term  has a computational cost equal to $4n_1n_2$ ops. The
algorithm proposed in \cite{ADNS-AR3}, that makes use of the
Anti-Reflective direct and inverse transforms, is less expensive
in the case of tests with few threshold values.\newline Hereafter,
the aim is to compare the truncated SVD with the truncated SD
restorations both in the case of Reflective and Anti-Reflective
BCs. Periodic BCs are not analyzed here, since Reflective and
Anti-Reflective BCs give better performances with respect to the
approximation of the image at the boundary.\newline In Table
\ref{tab:sigma2_cameraman_psfgauss} and
\ref{tab:sigma5_cameraman_psfgauss} we report the results obtained
by varying the dimension of the PSF support, the parameter $\rho$
related to the amount of the noise perturbation and the variance
of the considered Gaussian blur. As expected the optimal RRE
worses as the parameter $\rho$ increases and the Anti-Reflective
BCs show better performances in the case of  low noise levels. In
fact, for low $\rho$ values, the reduction of ringing artifacts is
significant, while the quality of the restoration for higher
$\rho$ values is essentially driven by the goal of noise
filtering.
%
% The results degrade with respect to  increasing noise levels.
%
Therefore, in such a case, the choice of the BCs becomes more an
more meaningless since it is not able to influence the image
restoration quality. Some examples of restored images are reported
in Figure \ref{fig:cameraman_psfgauss_d15s2ratio001-005}.\newline
%in Figures \ref{fig:cameraman_psfgauss_d15s2ratio001} and
%\ref{fig:cameraman_psfgauss_d15s2ratio005}.\newline
%-------------------------------------------------------------------------------------------------------------
\begin{table}
   \caption{Optimal RREs of truncated SVD and SD with reference to the true image in Figure
\ref{fig:cameraman_psfgauss_d15s2} (Gaussian blur
$\sigma_{i_1}=\sigma_{i_2}=2$).}
\label{tab:sigma2_cameraman_psfgauss}
\centering
\tiny%\footnotesize
\begin{tabular}{|l |llll|}
\hline \multicolumn{5}{|c|}{Reflective BCs}  \\\hline
PSF & \multicolumn{1}{|c}{5x5}      & \multicolumn{1}{c}{11x11}       &   \multicolumn{1}{c}{15x15}   &  \multicolumn{1}{c|}{21x21}\\
\hline \multicolumn{5}{|c|}{$\rho=$0} \\\hline
SVD  & 0.059164 &  0.087402  &     0.090742   &     0.093856  \\
%SD-s  &0.081363 &  0.087286 &  0.090746 &  0.093867 \\
SD  &0.043754 &  0.087400   &    0.090746 &  0.093867 \\
%---------------------------------------------------------
\hline \multicolumn{5}{|c|}{$\rho=$0.001} \\\hline
SVD  & 0.060278 & 0.091964  &  0.094468  & 0.097034 \\
SD   & 0.060278 & 0.091964  &  0.094476 &  0.097034\\
%---------------------------------------------------------
\hline \multicolumn{5}{|c|}{$\rho=$0.01} \\\hline
SVD  &0.091151 &  0.11214 &  0.11307  &  0.11495 \\
%SD-s  &0.098805 &  0.112140 &  0.113070 0&     0.11495 \\
SD  &0.091152 &  0.11214 &  0.11307 &  0.11495 \\
\hline \multicolumn{5}{|c|}{$\rho=$0.05} \\\hline
SVD  &0.11635 &  0.13356 &  0.13508  &  0.13739 \\
%SD-s  &0.116350 &  0.133560 &  0.135100  &     0.13739 \\
SD &0.11635 &  0.13356 &  0.13510  &     0.13739 \\
\hline \multicolumn{5}{|c|}{$\rho=$0.1} \\\hline
SVD  &0.13024 &  0.14607 &  0.14746 &  0.15047 \\
%SD-s  &0.130240 &  0.146070 &  0.147460 &  0.15047 \\
SD &0.13024 &  0.14607 &  0.14746 &  0.15047 \\
\hline
\end{tabular} \
\begin{tabular}{|l |llll|}
\hline \multicolumn{5}{|c|}{Anti-Reflective BCs}  \\\hline
PSF & \multicolumn{1}{|c}{5x5}      & \multicolumn{1}{c}{11x11}       &   \multicolumn{1}{c}{15x15}   &  \multicolumn{1}{c|}{21x21}\\
\hline \multicolumn{5}{|c|}{$\rho=$0} \\\hline
SVD  &0.039165 &  0.064081 &  0.086621 &  0.087237 \\
%SD-s  &0.079720 &  0.072589 &  0.083043 &  0.083521 \\
SD  &0.038316 &  0.063114 &  0.083043 &  0.083521 \\
%---------------------------------------------------------
\hline \multicolumn{5}{|c|}{$\rho=$0.001} \\\hline
SVD  & 0.062182 & 0.094237  & 0.098897  &  0.10042 \\
SD   & 0.059617 & 0.089105  & 0.092814  &  0.094343\\
%---------------------------------------------------------
\hline \multicolumn{5}{|c|}{$\rho=$0.01} \\\hline
SVD  &0.096049 &  0.12231 &  0.12403 &  0.12536 \\
%SD-s  &0.098947 &  0.112300  &     0.113430 &  0.11495 \\
SD  &0.091383 &  0.11230  &     0.11343 &  0.11495 \\
\hline \multicolumn{5}{|c|}{$\rho=$0.05} \\\hline
SVD  &0.12791 &  0.15070  &     0.15188 &  0.15492 \\
%SD-s  &0.116660 &  0.134140 &  0.135700  &     0.13816 \\
SD  &0.11666 &  0.13414 &  0.13570  &     0.13816 \\
\hline \multicolumn{5}{|c|}{$\rho=$0.1} \\\hline
SVD  &0.14399 &  0.16756 &  0.16964 &  0.17225 \\
%SD-s  &0.130830 &  0.147090 &  0.148520 &  0.15162 \\
SD &0.13083 &  0.14709 &  0.14852 &  0.15162 \\
\hline
\end{tabular}
\end{table}
%%%%%%%%%%%%%%%%%%%%%%%%%%%%%%%%%%%%%%%%%%%%%%%%%%%%%%%%%%%%%%%%%%%%%%%%%%%%%%%%%%%%%
%-------------------------------------------------------------------------------------------------------------
\begin{table}
   \caption{Optimal RREs of truncated SVD and SD with reference to the true image in Figure \ref{fig:cameraman_psfgauss_d15s2}
   (Gaussian blur  $\sigma_{i_1}=\sigma_{i_2}=5$).}
 \label{tab:sigma5_cameraman_psfgauss}
\centering
\tiny%\footnotesize
\begin{tabular}{|l |llll|}
\hline \multicolumn{5}{|c|}{Reflective BCs}  \\\hline
PSF & \multicolumn{1}{|c}{5x5}      & \multicolumn{1}{c}{11x11}       &   \multicolumn{1}{c}{15x15}   &  \multicolumn{1}{c|}{21x21}\\
\hline \multicolumn{5}{|c|}{$\rho=$0} \\\hline
SVD    & 0.063387 &    0.081274 &  0.097351 &  0.14634 \\
%SD-s     & 0.091528 &    0.129860 &  0.138430 &  0.154040 \\
SD   & 0.045365 &    0.081274 &  0.096387 &  0.14634 \\
%---------------------------------------------------------
\hline \multicolumn{5}{|c|}{$\rho=$0.001} \\\hline
SVD  & 0.063915 & 0.096243 & 0.11449 & 0.15217 \\
SD   & 0.063915 & 0.096274 & 0.11449 & 0.15217 \\
%---------------------------------------------------------
\hline \multicolumn{5}{|c|}{$\rho=$0.01} \\\hline
SVD    & 0.089032  &   0.13343 &   0.14947 &   0.17397 \\
%SD-s     & 0.106460  &   0.150990 &   0.162430 &   0.173970 \\
SD   & 0.089032  &   0.13343 &   0.14946 &   0.17397 \\
\hline \multicolumn{5}{|c|}{$\rho=$0.05} \\\hline
SVD    & 0.12203 & 0.16002 &   0.17339 &   0.18335 \\
%SD-s     & 0.121610 & 0.160950 &   0.173450 &   0.183350 \\
SD   & 0.12203 & 0.16002 &   0.17339 &   0.18335 \\
\hline \multicolumn{5}{|c|}{$\rho=$0.1} \\\hline
SVD    & 0.13412 & 0.16793 &   0.17963 &   0.19057 \\
%SD-s     & 0.134120 & 0.167930 &   0.179630 &   0.190570 \\
SD   & 0.13412 & 0.16793 &   0.17963 &   0.19057 \\
\hline
\end{tabular}
\
\begin{tabular}{|l |llll|}
\hline \multicolumn{5}{|c|}{Anti-Reflective BCs} \\\hline
PSF & \multicolumn{1}{|c}{5x5}      & \multicolumn{1}{c}{11x11}       &   \multicolumn{1}{c}{15x15}   &  \multicolumn{1}{c|}{21x21}\\
\hline \multicolumn{5}{|c|}{$\rho=$0} \\\hline
SVD  & 0.040214 &    0.079543  &  0.088224 &  0.13686 \\
%SD-s   & 0.089972 &    0.127050  &  0.139120 &  0.145180 \\
SD & 0.039437 &    0.078970  &  0.088832 &  0.13129 \\
%---------------------------------------------------------
\hline \multicolumn{5}{|c|}{$\rho=$0.001} \\\hline
SVD  & 0.068197 & 0.095808 & 0.11522  & 0.15767 \\
SD   & 0.063575 & 0.093247 & 0.1127   & 0.14893 \\
%---------------------------------------------------------
\hline \multicolumn{5}{|c|}{$\rho=$0.01} \\\hline
SVD  & 0.09412 &    0.14482 &   0.16825 &   0.21148 \\
%SD-s   & 0.106400 &    0.152260 &   0.163760 &   0.174460 \\
SD & 0.089038 &    0.13611 &   0.15270  &   0.17446 \\
\hline \multicolumn{5}{|c|}{$\rho=$0.05} \\\hline
SVD & 0.13553 & 0.18563 &   0.21006 &   0.22962 \\
%SD-s& 0.121820 & 0.162090 &   0.174200  &   0.184140 \\
SD & 0.12253 & 0.16269 &   0.17439 &   0.18414 \\
\hline \multicolumn{5}{|c|}{$\rho=$0.1} \\\hline
SVD & 0.15010  & 0.20164 &   0.22256 &   0.23960  \\
%SD-s& 0.134870 & 0.169160 &   0.180880 &   0.192180 \\
SD& 0.13487 & 0.16916 &   0.18088 &   0.19218\\
\hline
\end{tabular}
\end{table}
%%%%%%%%%%%%%%%%%%%%%%%%%%%%%%%%%%%%%%%%%%%%%%%%%%%%%%%%%%%%%%%%%%%
%
% COMMENTARE FIGURE!!!
%
%%%%%%%%%%%%%%%%%%%%%%%%%%%%%%%%%%%%%%%%%%%%%%%%%%%%%%%%%%%%%%%%%%%
%%%%%%%%%%%%%%%%%%%%%%%%%%%%%%%%%%%%%%%%%%%%%%%%%%%%%%%%%%%%%%%%%%%%%%%%%%%%%%%%%%%%%
%-------------------------------------------------------------------------------------------------------------
%%%%%%%%%%%%%%%%%%%%%%%%%%%%%%%%%%%%%%%%%%%%%%%%%%%%%%%%%%%%%%%%%%%%%%%%%%%%%%%%%%%%%
%
More impressive is the fact that SDs give better, or equal,
results with respect to those obtained by considering SVDs. This
numerical evidence is really interesting in the case of
Anti-Reflective BCs: despite the loss of the orthogonality
property in the spectral decomposition, the restoration results
are better than those obtained by considering SVD. Moreover, the
observed trend with respect to the Reflective BCs is also
conserved.\newline A further analysis refers to the so-called
Picard plots (see Figure
\ref{fig:cameraman_psfgauss_d15s2ratio001-005_picard_plot}),
%(see Figures
%\ref{fig:cameraman_psfgauss_d15s2ratio001_picard_plot} and
%\ref{fig:cameraman_psfgauss_d15s2ratio005_picard_plot}),
%
where the coefficients $|u_k^T g|$, or $|\tilde v_k g|$, (black
dots) are compared with the singular values $\sigma_k$, or the
absolute values of the eigenvalues $|\lambda_k|$, (red line). As
expected, initially these coefficients decrease faster than
$\sigma_k$, or $|\lambda_k|$, while afterwards  they level off at
a plateau determined by the level of the noise in the
image.\newline The threshold of this  change of behavior is in
good agreement with the optimal $k$ value obtained in the
numerical test by monitoring the RRE.\newline Moreover, notice
that the  Picard plots related to the SDs are quite in agreement
with those corresponding to SVDs. In the case of the
Anti-Reflective SD we observe an increasing data dispersion with
respect to the plateau, but the correspondence between the
threshold and the chosen optimal $k$ is still preserved.\newline
The computational relevance of this result is due to the
significant lower computational cost required by the
Anti-Reflective SDs with respect to the corresponding
SVDs.\newline
%%%%%%%%%%%%%%%%%%%%%%%%%%%%%%%%%%%%%%%%%%%%%%%%%%%%%%%%%%%%%%%%%%%%%%%%%%%%%%%%%%%%%
%%%%%%%%%%%%%%%%%%%%%%%%%%%%%%%%%%%%%%%%%%%%%%%%%%%%%%%%%%%%%%%%%%%%%%%%%%%%%%%%%%%%%
%-------------------------------------------------------------------------------------------------------------
\begin{table} %FF
\caption{Optimal RREs of truncated SDs with reference to the true
image in Figure \ref{fig:cameraman_psfgauss_d15s2} (Out-of-Focus
blur).}
\label{tab:cameraman_psfFF}
\centering
\tiny%\footnotesize
\begin{tabular}{|l |llll|llll|}
\hline \multicolumn{5}{|c|}{Reflective  BCs} &
\multicolumn{4}{|c|}{Anti-Reflective  BCs}
\\\hline
PSF & \multicolumn{1}{|c}{5x5}      & \multicolumn{1}{c}{11x11}
&   \multicolumn{1}{c}{15x15}   &  \multicolumn{1}{c|}{21x21}
& \multicolumn{1}{|c}{5x5}      & \multicolumn{1}{c}{11x11}       &   \multicolumn{1}{c}{15x15}   &  \multicolumn{1}{c|}{21x21} \\
\hline %\multicolumn{5}{|c|}{$\rho=$0} \\\hline
$\rho=$0& 0.072593 &    0.084604 &  0.088323 &  0.096479 & 0.072821 &    0.085366 &  0.091252 &  0.099293 \\
%---------------------------------------------------------
%\hline
%\multicolumn{5}{|c|}{$\rho=$0.001} \\\hline
$\rho=$0.001&   0.072671 & 0.085809 & 0.091035 & 0.10436 & 0.072904 & 0.086643 & 0.093929  & 0.10752  \\
%---------------------------------------------------------
%\hline
%\multicolumn{5}{|c|}{$\rho=$0.01} \\\hline
$\rho=$0.01& 0.080016 &    0.12255 &  0.13569 &  0.15276 & 0.080427 &    0.12316 &  0.13803 &  0.15683 \\
%\hline
%\multicolumn{5}{|c|}{$\rho=$0.05} \\\hline
$\rho=$0.05& 0.10645  &    0.15365 &  0.16810 &  0.18777 & 0.10685  &    0.15571 &  0.17147 &  0.19172 \\
%\hline
%\multicolumn{5}{|c|}{$\rho=$0.1} \\\hline
$\rho=$0.1& 0.12089  &    0.16314 &  0.17836 &  0.20471 & 0.12147  &    0.16482 &  0.17987 &  0.20829 \\
\hline
\end{tabular}
\end{table}
%%%%%%%%%%%%%%%%%%%%%%%%%%%%%%%%%%%%%%%%%%%%%%%%%%%%%%%%%%%%%%%%%%%%%%%%%%%%%%%%%%%%%
%
\noindent Lastly, Table  \ref{tab:cameraman_psfFF} reports the
spectral filtering results obtained in the case of Out-of-Focus
blur by varying the dimension of the PSF support and the parameter
$\rho$ related to the noise perturbation. The RRE follows the same
trend observed in the case of Gaussian blur. Other image
restoration tests with different gray-scale images have been
considered in \cite{FabrizioRossi}.\newline A more interesting
remark again pertains  the computational cost. Since the
Out-of-Focus PSF is not separable, but the transforms are, the use
of SDs related to Reflective or Anti-Reflective BCs allows to
exploit the tensorial nature of the corresponding transforms, both
with respect to the computation of the eigenvalues and of the
eigenvectors (or of the Reflective and Anti-Reflective
transforms). \newline
%-----------------------------------sub-section----------------------------------------------------------------
\subsubsection{COLOR IMAGES IN THE CASE OF CROSS-CHANNEL BLURRING} \label{sez:cross_channel_blurring_numerical}
\begin{figure}%[h!]
   \centering
   \includegraphics[width=\textwidth]{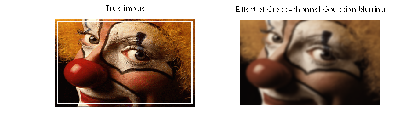}
   \caption{True image (FOV is delimited by a white frame) and cross-channel blurred image
   with the Gaussian PSF with support
$15\times 15$ and $\sigma_{i_1}=\sigma_{i_2}=2$ and matrix
$A_{\mathrm{color}}$ in (\ref{eq:def_Acolor}).}
   \label{fig:clown_psfgauss_d15s2}
\end{figure}

%-------------------------------------------------------------------------------------------------------------
Here, we analyze some restoration tests in the case the template
color image  reported in Figure \ref{fig:clown_psfgauss_d15s2}, by
assuming the presence of a cross-channel blurring phenomenon
modelled according to (\ref{eq:sistema_quadrato_color}). The
entity of this mixing effect is chosen according to the matrix
\begin{equation} \label{eq:def_Acolor}
A_{\mathrm{color}}= \left [
\begin{array}{ccc}
0.7& 0.2& 0.1\\
0.25& 0.5& 0.25\\
0.15& 0.1& 0.75
\end{array}
\right ].
\end{equation}
%%%%%%%%%%%%%%%%%%%%%%%%%%%%%%%%%%%%%%%%%%%%%%%%%%%%%%%%%%%%%%%%%%%%%%%%%%%%%%%%%%%%%
%%%%%%%%%%%%%%%%%%%%%%%%%%%%%%%%%%%%%%%%%%%%%%%%%%%%%%%%%%%%%%%%%%%%%%%%%%%%%%%%%%%%%
%
In Figure \ref{fig:clown_psfgauss_d15s2} is also reported the
cross-channel blurred image with Gaussian PSF with support
$15\times 15$ and $\sigma_{i_1}=\sigma_{i_2}=2$. Notice that the
entity of the cross-channel blurring is not negligible, since the
whole image results to be darkened and the color intensities of
the additive RGB system are substantially altered.\newline Table
\ref{tab:clown_psfgauss_sigma2} reports the optimal RREs of
truncated SVDs and SDs obtained by varying the dimension of the
Gaussian PSF support and the parameter $\rho$ related to the
amount of the noise perturbation. It is worth stressing that we
tune the filtering strategy with respect the spectral information
given just by the matrix $A_n$, i.e., for any fixed $\sigma_k$ (or
$\lambda_k$) we simultaneously sum, or discard, the three
contribution on $f$ related to the three singular values of
$A_{\mathrm{color}}$. In fact, the magnitude of singular values of
the considered matrix $A_{\mathrm{color}}$ does not differ enough
to dramatically change the filtering information given just by
$A_n$. Nevertheless, also the comparison with the restoration
results obtained by considering a global ordering justifies this
approach.\newline
%
%Quite no new remark can be done with respect to the gray-scale
%case:
The color case behaves as the gray-scale one: as expected the
optimal RRE becomes worse as the parameter $\rho$ increases and
the Anti-Reflective SD shows better performances in the case of
low noise levels.
\newline
%%%%%%%%%%%%%%%%%%%%%%%%%%%%%%%%%%%%%%%%%%%%%%%%%%%%%%%%%%%%%%%%%%%%%%%%%%%%%%%%%%%%%
\begin{figure}%[bh!]
   \centering
   $\rho=0.01$ \\
   \includegraphics[width=\textwidth]{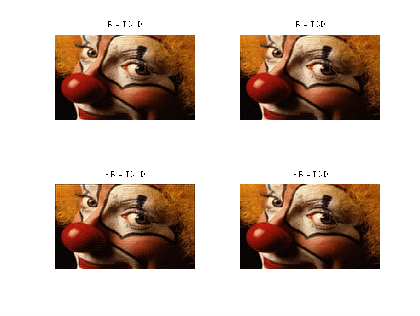}
%
%   \caption{Optimal restorations of truncated SVD and SD in the
%   case of Reflective and Anti-Reflective BCs with reference to Figure \ref{fig:clown_psfgauss_d15s2} - $\rho=0.01$.}
%   \label{fig:clown_psfgauss_d15s2ratio001}
%\end{figure}
%%%%%%%%%%%%%%%%%%%%%%%%%%%%%%%%%%%%%%%%%%%%%%%%%%%%%%%%%%%%%%%%%%%%%%%%%%%%%%%%%%%%%
%\begin{figure}
      \centering
   $\rho=0.05$ \\
   \includegraphics[width=\textwidth]{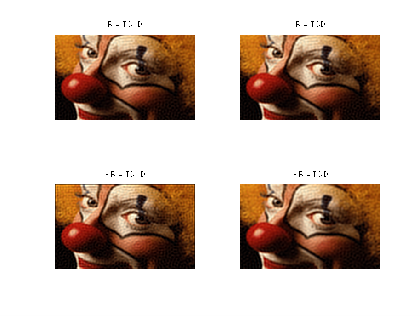}
   \caption{Optimal restorations of truncated SVD and SD in the
   case of Reflective and Anti-Reflective BCs with reference to Figure \ref{fig:clown_psfgauss_d15s2}
   (Cross-channel and Gaussian Blur $\sigma_{i_1}=\sigma_{i_2}=2$). % - $\rho=0.05$.
   }
   \label{fig:clown_psfgauss_d15s2ratio001-005}
%   \label{fig:clown_psfgauss_d15s2ratio005}
\end{figure}
%%%%%%%%%%%%%%%%%%%%%%%%%%%%%%%%%%%%%%%%%%%%%%%%%%%%%%%%%%%%%%%%%%%%%%%%%%%%%%%%%%%%%
In addition, by referring to Figure
\ref{fig:clown_psfgauss_d15s2ratio001-005},
% and \ref{fig:clown_psfgauss_d15s2ratio005},
we note that the truncated SVD in the case of Anti-Reflective BCs
shows a little more 'freckles' than the corresponding truncated
SVD in the case of Reflective BCs. Nevertheless, for low noise
levels, is just the Anti-Reflective SD that exhibits less
'freckles' than the Reflective SD.\newline
%%%%%%%%%%%%%%%%%%%%%%%%%%%%%%%%%%%%%%%%%%%%%%%%%%%%%%%%%%%%%%%%%%%%%%%%%%%%%%%%%%%%%%
\begin{table}%[t!]
\caption{Optimal RREs of truncated SVD and SD with reference to
the true image in Figure \ref{fig:clown_psfgauss_d15s2}
(Cross-channel and Gaussian Blur $\sigma_{i_1}=\sigma_{i_2}=2$).}
\label{tab:clown_psfgauss_sigma2}
\centering
\tiny%\footnotesize
\begin{tabular}{|l |llll|}
\hline \multicolumn{5}{|c|}{Reflective BCs}  \\\hline
PSF & \multicolumn{1}{|c}{5x5}      & \multicolumn{1}{c}{11x11}       &   \multicolumn{1}{c}{15x15}   &  \multicolumn{1}{c|}{21x21}\\
\hline \multicolumn{5}{|c|}{$\rho=$0} \\\hline
SVD &0.078276 &0.12114   &0.11654   &0.1178 \\
SD &0.078276  &0.12114   &0.11654   &0.1178 \\
%
%---------------------------------------------------------
\hline \multicolumn{5}{|c|}{$\rho=$0.001} \\\hline
SVD  & 0.078992 & 0.1212  & 0.11663 & 0.11792 \\
SD   & 0.078992 & 0.12119 & 0.11663 & 0.11792 \\
%---------------------------------------------------------
\hline \multicolumn{5}{|c|}{$\rho=$0.01} \\\hline
SVD &0.10152  &0.12396   &0.12088   &0.12198 \\
SD &0.10152   &0.12396   &0.12088   &0.12198 \\
\hline \multicolumn{5}{|c|}{$\rho=$0.05} \\\hline
SVD &0.12102   &0.13853   &0.13743   &0.13844 \\
SD &0.12102    &0.13853   &0.13743   &0.13844 \\
\hline \multicolumn{5}{|c|}{$\rho=$0.1} \\\hline
SVD &0.13437   &0.14898   &0.14854   &0.14947 \\
SD  &0.13437   &0.14898   &0.14854   &0.14947 \\
\hline
\end{tabular}
\begin{tabular}{|l |llll|}
\hline \multicolumn{5}{|c|}{Anti-Reflective BCs}  \\\hline
PSF & \multicolumn{1}{|c}{5x5}      & \multicolumn{1}{c}{11x11}       &   \multicolumn{1}{c}{15x15}   &  \multicolumn{1}{c|}{21x21}\\
\hline \multicolumn{5}{|c|}{$\rho=$0} \\\hline
SVD &0.076646  &0.1006   &0.1098   &0.10646 \\
SD  &0.074953  &0.098474 &0.10508  &0.10216 \\
%
%---------------------------------------------------------
\hline \multicolumn{5}{|c|}{$\rho=$0.001} \\\hline
SVD  &0.077394 & 0.10639 & 0.1111  & 0.11002 \\
SD   &0.075727 & 0.10233 & 0.10612 & 0.10443 \\
%---------------------------------------------------------
\hline \multicolumn{5}{|c|}{$\rho=$0.01} \\\hline
SVD &0.10431  &0.12695   &0.12624   &0.12779 \\
SD &0.10087   &0.11737   &0.11805   &0.118 \\
\hline \multicolumn{5}{|c|}{$\rho=$0.05} \\\hline
SVD &0.13017   &0.15075   &0.15063   &0.15166 \\
SD &0.12127    &0.13699   &0.13756   &0.13795 \\
\hline \multicolumn{5}{|c|}{$\rho=$0.1} \\\hline
SVD &0.1456   &0.16516   &0.16626   &0.16647 \\
SD &0.13507    &0.14796   &0.14955   &0.15018 \\
\hline
\end{tabular}
\end{table}
\subsection{TIKHONOV AND RE-BLURRING REGULARIZATIONS}
\begin{figure}
   \centering
   $\rho=0.01$ \\
   \vskip -1.5cm
   \includegraphics[width=\textwidth]{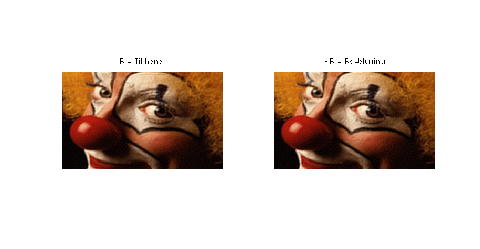}
%
%   \caption{Optimal RREs  of Tikhonov and Re-blurring methods and
%corresponding $\mu_{ott}$ with reference to the true image in
%Figure \ref{fig:clown_psfgauss_d15s2} - Cross-channel Gaussian
%blur $\sigma_{i_1}=\sigma_{i_2}=2$ - $\rho=0.01$.}
%   \label{fig:clown_psfgauss_tik_d15s2ratio001}
%\end{figure}
%%%%%%%%%%%%%%%%%%%%%%%%%%%%%%%%%%%%%%%%%%%%%%%%%%%%%%%%%%%%%%%%%%%%%%%%%%%%%%%%%%%%%%
%%%%%%%%%%%%%%%%%%%%%%%%%%%%%%%%%%%%%%%%%%%%%%%%%%%%%%%%%%%%%%%%%%%%%%%%%%%%%%%%%%%%%
%\begin{figure}
   \vskip -1.5cm
   \centering
   $\rho=0.05$ \\
   \vskip -1.5cm
   \includegraphics[width=\textwidth]{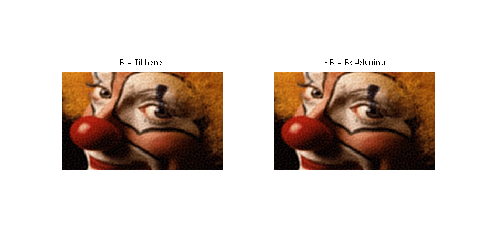}
   \vskip -1cm
   \caption{Optimal RREs  of Tikhonov and Re-blurring methods with reference to the true image in
Figure \ref{fig:clown_psfgauss_d15s2} (Cross-channel and Gaussian
blur $\sigma_{i_1}=\sigma_{i_2}=2$ - $\rho=0.05$).}
   \label{fig:clown_psfgauss_tik_d15s2ratio001-005}
\end{figure}

%-------------------------------------------------------------------------------------------------------------
\begin{table} %
\caption{Optimal RREs of Tikhonov and Re-blurring methods and
corresponding $\mu_{ott}$ with reference to the true image in
Figure \ref{fig:cameraman_psfgauss_d15s2} (Gaussian Blur
$\sigma_{i_1}=\sigma_{i_2}=2$).}
\label{tab:cameraman_psfGauss_Tikhnov_sigma2}
\centering
\tiny%\footnotesize
\begin{tabular}{|l |ll|ll|ll|ll|}
%\hline \multicolumn{5}{|c|}{Reflective  BCs}\\
\hline
PSF & \multicolumn{2}{|c|}{5x5}      & \multicolumn{2}{|c|}{11x11}       &   \multicolumn{2}{|c|}{15x15}   &  \multicolumn{2}{|c|}{21x21}\\
\hline \multicolumn{9}{|c|}{$\rho=$0} \\\hline
R  & 0.041015 & 4.1e-005    & 0.079044 & 9e-006 & 0.086386 & 1.1e-005   & 0.089556 & 1.6e-005\\
AR & 0.034237 & 1.1e-005    & 0.059465 & 1e-006 & 0.078963 & 1e-006 & 0.079805 & 1e-006 \\
%------------------------------------------------------------------------------------------------
\hline \multicolumn{9}{|c|}{$\rho=$0.001} \\\hline
R  &  0.050155 & 0.000188 &  0.087482 & 5.7e-005    &  0.090825 & 4.3e-005     &  0.093071 & 4.9e-005 \\
AR &  0.048556 & 0.000163 &  0.085279 & 4.6e-005    &  0.089388 & 3.3e-005     &  0.090821 & 3.3e-005 \\
%------------------------------------------------------------------------------------------------
\hline \multicolumn{9}{|c|}{$\rho=$0.01} \\\hline
R  & 0.083456 & 0.005555    & 0.10748  & 0.001786   & 0.10863 & 0.001678    & 0.11023  & 0.001573\\
AR & 0.083436 & 0.005536    & 0.10744  & 0.001792   & 0.10868 & 0.001691    & 0.11019  & 0.001575\\
\hline \multicolumn{9}{|c|}{$\rho=$0.05} \\\hline
R  & 0.12024  & 0.038152    & 0.12982  & 0.01929    & 0.13071 & 0.018417    & 0.13307  & 0.017892\\
AR & 0.12049  & 0.038379    & 0.13006  & 0.01957    & 0.13096 & 0.018669    & 0.1333   & 0.018105\\
\hline \multicolumn{9}{|c|}{$\rho=$0.1} \\\hline
R  & 0.14767  & 0.06587 & 0.14721  & 0.039231   & 0.14822 & 0.038181    & 0.15097  & 0.037893\\
AR & 0.14813  & 0.066251    & 0.14766  & 0.039707   & 0.14866 & 0.038644    & 0.15144  & 0.038296\\
\hline
\end{tabular}
\end{table}
\begin{table} %
\caption{Optimal RREs  of Tikhonov and Re-blurring methods and
corresponding $\mu_{ott}$ with reference to the true image in
Figure \ref{fig:cameraman_psfgauss_d15s2} (Out-of-Focus blur).}
\label{tab:cameraman_psfFF_Tikhnov}
\centering
\tiny%\footnotesize
\begin{tabular}{|l |ll|ll|ll|ll|}
%\hline \multicolumn{5}{|c|}{Reflective  BCs}\\
\hline
PSF & \multicolumn{2}{|c|}{5x5}      & \multicolumn{2}{|c|}{11x11}       &   \multicolumn{2}{|c|}{15x15}   &  \multicolumn{2}{|c|}{21x21}\\
\hline \multicolumn{9}{|c|}{$\rho=$0} \\\hline
R & 0.031422 & 0.000172 & 0.05346 & 6.9e-005    & 0.060954 & 3.5e-005   & 0.074785 & 2.7e-005\\
AR & 0.036213 & 0.000302    & 0.051236 & 6.8e-005   & 0.06683 & 5.7e-005    & 0.084482 & 5.8e-005 \\
%------------------------------------------------------------------------------------------------
\hline \multicolumn{9}{|c|}{$\rho=$0.001} \\\hline
R  &  0.034441 & 0.000271  &  0.061465 & 0.000145 &  0.073751 & 0.000101   & 0.09074  & 7.9e-005\\
AR &  0.038313 & 0.000402  &  0.059957 & 0.000138 &  0.076695 & 0.000126   & 0.095274 & 0.000106 \\
%------------------------------------------------------------------------------------------------
\hline \multicolumn{9}{|c|}{$\rho=$0.01} \\\hline
R & 0.069647 & 0.008493 & 0.11361 & 0.004117    & 0.12881 & 0.003037    & 0.14914 & 0.001873\\
AR & 0.070384 & 0.008923    & 0.11404 & 0.00422 & 0.12982 & 0.003139    & 0.15061 & 0.001969\\
\hline \multicolumn{9}{|c|}{$\rho=$0.05} \\\hline
R & 0.12204 & 0.053687  & 0.1532 & 0.030719 & 0.16614 & 0.022121    & 0.18769 & 0.01346\\
AR & 0.12256 & 0.05423  & 0.15402 & 0.031574    & 0.16739 & 0.023213    & 0.18933 & 0.014472\\
\hline \multicolumn{9}{|c|}{$\rho=$0.1} \\\hline
R & 0.16366 & 0.092379  & 0.17357 & 0.055919    & 0.1829 & 0.042944 & 0.20323 & 0.028803\\
AR & 0.16433 & 0.093069 & 0.17485 & 0.057326    & 0.18457 & 0.044901    & 0.20511 & 0.031011\\
\hline
\end{tabular}
\end{table}
%
%-------------------------------------------------------------------------------------------------------------
By considering a Gaussian blurring of the true image reported in
Figure \ref{fig:cameraman_psfgauss_d15s2}, Table
\ref{tab:cameraman_psfGauss_Tikhnov_sigma2} compares the optimal
RRE obtained in the case of the Tikhonov method for Reflective BCs
and of the Re-blurring method for Anti-Reflective BCs. In
addition, in Table \ref{tab:cameraman_psfFF_Tikhnov}, the same
comparison refers to the case of the Out-of-Focus PSF.\newline As
expected, the RRE deteriorates as the dimension of the noise level
or the dimension of the PSF support increases. Notice also that
the gap between the Reflective and Anti-Reflective BCs is reduced
also for low noise levels. Further numerical tests can be found in
\cite{DEMS-IP-2006,ADS-AR1}.\newline Lastly, we focus our
attention on the case of the color image in Figure
\ref{fig:clown_psfgauss_d15s2}. The image restorations have been
obtained by considering the transformation procedure outlined at
the end of Section \ref{sez:cross_channel_blurring}. Despite the
RREs in Table
\ref{tab:cameraman_psfGausssigma2_Tikhnov_cross_channel} are
bigger than in the gray-scale case, the perception of the image
restoration quality is very satisfying and a little less
'freckles' than in the corresponding SDs and SVDs are observed
(see Figure \ref{fig:clown_psfgauss_tik_d15s2ratio001-005}).
Notice, also that the lack of orthogonality in the $S_{3n}$
transform related to the Anti-reflective BCs does not deteriorate
the performances of the restoration.
\begin{table} %
\caption{Optimal RREs  of Tikhonov and Re-blurring methods and
corresponding $\mu_{ott}$ with reference to the true image in
Figure \ref{fig:clown_psfgauss_d15s2} (Cross-channel and Gaussian
Blur $\sigma_{i_1}=\sigma_{i_2}=2$).}
\label{tab:cameraman_psfGausssigma2_Tikhnov_cross_channel}
\centering
\tiny%\footnotesize
\begin{tabular}{|l |ll|ll|ll|ll|}
%\hline \multicolumn{5}{|c|}{Reflective  BCs}\\
\hline
PSF & \multicolumn{2}{|c|}{5x5}      & \multicolumn{2}{|c|}{11x11}       &   \multicolumn{2}{|c|}{15x15}   &  \multicolumn{2}{|c|}{21x21}\\
\hline \multicolumn{9}{|c|}{$\rho=$0} \\\hline
R & 0.069148 & 0.000203   & 0.11508 & 0.001204    & 0.1123 & 0.000717 & 0.11335 & 0.000726 \\
AR & 0.062854 & 0.000102   & 0.091232 & 7e-006 & 0.1014 & 4.4e-005 & 0.098266 & 1.5e-005 \\
%------------------------------------------------------------------------------------------------
\hline \multicolumn{9}{|c|}{$\rho=$0.001} \\\hline
R  &  0.071259 & 0.000312 &  0.11515  & 0.001228 &  0.11239 & 0.000744 &   0.11347 & 0.000755 \\
AR &  0.066734 & 0.000209 &  0.098658 & 5.8e-005 &  0.10276 & 7.7e-005 &   0.10111 & 4.5e-005 \\
%------------------------------------------------------------------------------------------------
\hline \multicolumn{9}{|c|}{$\rho=$0.01} \\\hline
R & 0.094871 & 0.004975   & 0.11896 & 0.002919    & 0.11712 & 0.002421    & 0.1182  & 0.002459 \\
AR & 0.094458 & 0.004841   & 0.1144 & 0.001884 & 0.11481 & 0.00184 & 0.11507 & 0.001755 \\
\hline \multicolumn{9}{|c|}{$\rho=$0.05} \\\hline
R & 0.13209 & 0.029798    & 0.13662 & 0.015305    & 0.13599 & 0.014896    & 0.13669 & 0.014824 \\
AR & 0.13239 & 0.029944    & 0.13561 & 0.014992    & 0.13593 & 0.014772    & 0.13611 & 0.014595 \\
\hline \multicolumn{9}{|c|}{$\rho=$0.1} \\\hline
R & 0.16281 & 0.051315    & 0.15543 & 0.029068    & 0.15547 & 0.028822    & 0.15586 & 0.02868 \\
AR & 0.16341 & 0.051659    & 0.15526 & 0.029213    & 0.15602 & 0.029029    & 0.15588 & 0.02872 \\
\hline
\end{tabular}
\end{table}
%
%%%%%%%%%%%%%%%%%%%%%%%%%%%%%%%%%%%%%%%%%%%%%%%%%%%%%%%%%%%%%%%%%%%%%%%%%%%%%%%%%%%%%
%%%%%%%%%%%%%%%%%%%%%%%%%%%%%%%%%%%%%%%%%%%%%%%%%%%%%%%%%%%%%%%%%%%%%%%%%%%%%%%%%%%%%
%-----------------------------------Section----------------------------------------------------------------
\section{CONCLUSIONS} \label{sez:conclusioni}
In this paper we have analyzed and compared SD and SVD filtering
methods in the case both of Reflective and Anti-Reflective BCs.
Numerical evidence is given of the good performances achievable
through SDs and with a substantially lower computational cost with
respect to SVDs. In addition, the tensorial structure of the
Reflective and Anti-Reflective SDs can be exploited in depth also
in the case of not separable PSFs. \newline A special mention has
to be done to the fact that the loss of orthogonality of the
Anti-Reflective transform does not seems to have any consequence
on the trend of the image restoration results. The analysis in the
case of cross-channel blurring in color images allows to confirm
the quoted considerations. Finally, the Re-blurring regularizing
method has been re-interpreted as a standard Tikhonov
regularization method in the space transformed by means of
$\widetilde T_n$. Some numerical tests highlight the image
restoration performances, also in the case of cross-channel
blurring. \newline Future works will concern the analysis of
effective strategies allowing to properly choose the optimal
regularizing parameters in the Anti-Reflective BCs case.\newline
%-----------------------------------Section----------------------------------------------------------------

\end{document}